\documentclass[12pt]{amsart}%
\usepackage{amsmath,amscd}%
\usepackage{amsfonts}%
\usepackage{amssymb}%
\usepackage{graphicx}
\usepackage{epsfig}

\newtheorem{theorem}{Theorem}
\theoremstyle{plain}

\newtheorem{corollary}{Corollary}

\newtheorem{definition}{Definition}

\newtheorem{lemma}{Lemma}

\newtheorem{proposition}[theorem]{Proposition}
\newtheorem{remark}{Remark}

\numberwithin{equation}{section}

\begin{document}

\title{Moduli Spaces and Multiple Polylogarithm Motives}
\author{Qingxue Wang}
\maketitle

\begin{abstract}
In this paper, we give a natural construction of mixed Tate
motives whose periods are a class of iterated integrals which
include the multiple polylogarithm functions. Given such an
iterated integral, we construct two divisors $A$ and $B$ in the
moduli spaces $\mathcal{\overline{M}}_{0,n}$ of $n$-pointed stable
curves of genus $0$, and prove that the cohomology of the pair
$(\mathcal{\overline{M}}_{0,n}-A,B-B\cap A)$ is a framed mixed
Tate motive whose period is that integral. It generalizes the
results of A. Goncharov and Yu. Manin for multiple $\zeta$-values.
Then we apply our construction to the dilogarithm and calculate
the period matrix which turns out to be same with the canonical
one of Deligne.
\end{abstract}

\section{Introduction}
\subsection{Multiple Polylogarithms}
The multiple polylogarithm functions were defined in Goncharov's
paper \cite{Gon0} as the following power series:
\begin{equation}\label{L}
Li_{n_1, \dots ,n_m}(x_1,\dots,x_m)=\sum_{0<k_1<k_2< \dots
<k_m}\frac{x_1^{k_1}x_2^{k_2}\dots
x_m^{k_m}}{k_1^{n_1}k_2^{n_2}\dots k_{m}^{n_{m}}}
\end{equation}
where the $x_i$ are in the unit disk of the complex plane for
$i=1,\dots,m$ and $n_1\geq 1,\dots,n_{m-1}\geq 1, n_m\geq 2$ are
positive integers. For $m=1$, we get the classical $n$-th
polylogarithm which was first introduced by Leibniz \cite{L} in
1696:
\begin{equation}\label{Pl}
Li_n(z)=\sum_{k=1}^{+\infty}\frac{z^k}{k^n}\ , \quad |z|\leq 1
\end{equation}
And for $x_1= \dots =x_m=1$, we obtain the multiple $\zeta$-values
which were first studied by Euler \cite{Eul}:
\begin{equation}\label{mz}
\zeta(n_1, \dots ,n_m)=\sum_{0<k_1<k_2< \dots
<k_m}\frac{1}{k_1^{n_1}k_2^{n_2}\dots k_{m}^{n_{m}}}\ .\\
\end{equation}
Moveover, multiple polylogarithms can be represented as iterated
integrals. Recall that iterated integrals are defined as follows.
Let $\omega_1,\dots,\omega_n$ be smooth one-forms on a manifold
$M$ and $\gamma:[0,1] \rightarrow M$ be a piecewise smooth path.
Then we define inductively as follows:
\[ \int_{\gamma}{\omega_1 \circ \dots \circ
\omega_n}:=\int_0^1{(\int_{\gamma_t}{\omega_1 \circ \dots \circ
\omega_{n-1}})\gamma^{*}\omega_n}
\]
where $\gamma_t$ is the restriction of $\gamma$ on $[0,t]$ and
$\int_{\gamma_t}{\omega_1 \circ \dots \circ \omega_{n-1}}$ is a
function of $t$ on $[0,1]$. More explicitly, it can be computed in
the following way:
\[\int_{\gamma}{\omega_1 \circ \dots \circ \omega_n}=\int_{0\leq
t_1\leq \dots \leq t_n\leq 1}{f_1(t_1)\, dt_1\wedge \dots \wedge
f_n(t_n)\, dt_n}
\]
where $f_i(t)dt=\gamma^{*}\omega_i$ are the pullback one-forms on
$[0,1]$, $i=1,\dots ,n$. For example,
\[
\zeta(2)=\sum_{n=1}^{+\infty}{\frac{1}{n^2}}=\int_0^1{\frac{dt}{1-t}
\circ \frac{dt}{t}}=\int_{0\leq t_1\leq t_2\leq
1}{\frac{dt_1}{1-t_1}\wedge \frac{dt_2}{t_2}}
\]In \cite[Chap2]{Gon1},
the following formula was proved:
\begin{equation*}
Li_{n_1, \dots ,n_m}(x_1, \dots ,x_m)=
\end{equation*}
\begin{equation}\label{mul-li}
(-1)^m \; \int_{0}^{1}{\underbrace{\frac{dt}{t-(x_{1} \dots
x_{m})^{-1}} \circ \frac{dt}{t} \circ \dots \circ
\frac{dt}{t}}_{n_{1}\quad times} \circ  \dots  \circ
\underbrace{\frac{dt}{t-x_{m}^{-1}}\circ \frac{dt}{t} \circ \dots
\circ \frac{dt}{t},}_{n_{m}\quad times}}
\end{equation}
This formula also provides the analytic continuation of multiple
polylogarithms.
\subsection{Moduli Spaces}
We denote by $\mathcal{\overline{M}}_{0,S}$ the moduli space of
$S$-labeled pointed stable curves of genus $0$, where $S$ is a
finite set. It's been studied by Grothendieck \cite{De4}, Deligne,
Mumford, Knudsen \cite{Kn} and many others. It is defined over
$\mathbb{Z}$. Roughly speaking, a complex point of
$\mathcal{\overline{M}}_{0,S}$ is a tree of complex projective
lines with $|S|$ distinct smooth points marked by the set $S$.
Here $|S|$ denotes the cardinality of the set $S$. We know that it
is a smooth irreducible projective variety of complex dimension
$|S|-3$. Moreover, $\mathcal{\overline{M}}_{0,S}(\mathbb{C})$
provides a natural compactification of the space
$\mathcal{M}_{0,S}(\mathbb{C})$ of $|S|$ distinct points on
$\mathbb{CP}^1$ modulo automorphisms of $\mathbb{CP}^1$. By
cross-ratio, $\mathcal{M}_{0,S}(\mathbb{C})$ is isomorphic to
\[\{(x_1,\dots,x_n)\in(\mathbb{CP}^1)^n|x_i\ne
x_j,i\ne j\; ;x_k\ne 0,1,\infty\,\  k=1,\dots,n\}, \; n=|S|-3.\]
The boundary $\partial
\mathcal{\overline{M}}_{0,S}:=\mathcal{\overline{M}}_{0,S}-\mathcal{M}_{0,S}$
is a normal crossing divisor. It can be described by the
combinatorial data of the set $S$. For more detail, see section
$2$. Now given any subset $S_1\subset S$, with $|S_1|\geq 3$,
there is a contraction morphism:
\[\pi_{S_1}: \mathcal{\overline{M}}_{0,S}\rightarrow
\mathcal{\overline{M}}_{0,S_1}\] which contracts stably all
sections but those marked by $S_1$. In particular, for any $s_0\in
S$, let $S^{\prime}=S\setminus\{s_0\}$, the contraction morphism
\[\pi:\mathcal{\overline{M}}_{0,S}\rightarrow
\mathcal{\overline{M}}_{0,S^{\prime}}\] is the universal
$S^{\prime}$-labeled curve with universal sections $\sigma_i$, for
each $i\in S^{\prime}$. For more information and proofs of
$\mathcal{\overline{M}}_{0,S}$, we refer to \cite{Kn,Ke,Man2}.

If the set $S=\{1,2,\dots,n\}$, we'll denote this space by
$\mathcal{\overline{M}}_{0,n}$. And $\mathcal{\overline{M}}_{0,S}$
is non-canonically isomorphic to $\mathcal{\overline{M}}_{0,|S|}$.
From now on, we'll take $S=\{0,s_1,s_2,\dots,$\\
$s_n,1,\infty\}$ and fix the cyclic order $\rho:
0<s_1<s_2<\dots<s_n<1<\infty<0$ on $S$ unless otherwise stated.
\subsection{Main results}
In \cite{GM}, for each multiple $\zeta$-value (\ref{mz}), the
authors construct two divisors $A$ and $B$ of
$\mathcal{\overline{M}}_{0,S}$ and then show that
\[H^{n}(\mathcal{\overline{M}}_{0,S}-A, B-B \cap A)\]
is a framed mixed Tate motive whose period is this value. In the
end of that paper, they suggest to generalize their results to the
following convergent iterated integral:
\begin{equation}\label{int1}
I_{\gamma}(a_1, \dots , a_n):= \int_{\gamma}{\frac{dt}{t-a_1}\circ
\dots \circ\frac{dt}{t-a_{n}}} \quad a_{1}\neq0,a_{n}\neq1
\end{equation}
where $\gamma:[0,1] \rightarrow \mathbb{C}$ is a piecewise smooth
simple path from $0$ to $1$ and $a_i\notin \gamma((0,1))$,
$i=1,\dots,n$. In particular, by the formula (\ref{mul-li}),
multiple polylogarithms are of this type. In this paper, we show
that the analogous results hold for the iterated integral
(\ref{int1}).

In section $2$, we review the basic combinatorial facts about the
boundary divisors of $\mathcal{\overline{M}}_{0,S}$ and the stable
$2$-partitions of the set $S$. Next, we briefly recall the divisor
$B_n$ in $\mathcal{\overline{M}}_{0,S}$ which was introduced in
\cite{GM}, and then prove some interesting combinatorial
properties of $B_n$. In the end, we proceed to study in detail
some non-boundary divisors of $\mathcal{\overline{M}}_{0,S}$ which
we'll use later on.

In section $3$, for the integral (\ref{int1}), we define a
meromorphic differential form $\Omega_S(\vec{a})$ of
$\mathcal{\overline{M}}_{0,S}(\mathbb{C})$. Let $A_S(\vec{a})$ be
its divisor of singularities in
$\mathcal{\overline{M}}_{0,S}(\mathbb{C})$. We explicitly
determine the divisor $A_S(\vec{a})$. Then we use it to prove the
key proposition that \emph{the divisor $A_S(\vec{a})$ does not
contain any $k$-dimensional face of the divisor $B_n$, $0\leq
k\leq n$}.

In section $4$, we review the definitions of framed Hodge-Tate
structure and its period, and discuss their basic properties.
Finally combining with all the information of $A_S(\vec{a})$ and
$B_n$ in section $3$ and $4$, we can prove:
\begin{theorem}\label{main}
Let $\vec{a}=(a_1,\dots,a_n)$. For the iterated integral
$I_{\gamma}(a_1,\dots,a_n)$ of (\ref{int1}), $a_i\in
\mathbb{C},a_{1}\neq0,a_{n}\neq1$, there exists two divisors
$A_S(\vec{a})$ and $B_n$ in $\mathcal{\overline{M}}_{0,S}$,
$|S|=n+3$, such that
\[
H^{n}(\mathcal{\overline{M}}_{0,S}-A_S(\vec{a}), B_{n}-B_{n} \cap
A_S(\vec{a}))    \qquad (**)\] carries an $n$-framed
Hodge-Tatestructure with two canonical frames
\begin{equation*}
{[\Omega_S(\vec{a})]}\in
Gr_{2n}^{W}H^{n}(\mathcal{\overline{M}}_{0,S}-A_S(\vec{a})) ;\;
{[\Delta_{B}(\gamma)]}\in
(Gr_{0}^{W}H^{n}(\mathcal{\overline{M}}_{0,S},B_n))^{\vee}.
\end{equation*}
and the period with respect to these frames is exactly the
iterated integral $I_{\gamma}(a_1,\dots,a_n)$, where
$[\Omega_S(\vec{a})]$ is a meromorphic $n$-form on
$\mathcal{\overline{M}}_{0,S }$, and $\Delta_{B}(\gamma)$ is a
relative $n$-cycle.

Furthermore, if the $a_i$ are elements of a number field $F$,
$i=1,\dots,n$, then $(**)$ is a framed mixed Tate motive over $F$.
\end{theorem}

In section $5$, we apply our construction to the dilogarithm.
Namely, we consider the following integral:
\[Li_{2}(z)=-\int_{0\leq t_{1}\leq t_{2}\leq
1}{\frac{dt_{1}}{t_{1}-z^{-1}}\wedge \frac{dt_{2}}{t_{2}}}\] By
the Theorem \ref{main} above,
$H^{2}(\mathcal{\overline{M}}_{0,5}-A(\vec{a}), B_2-B_2 \cap
A(\vec{a}))$ carries a $2$-framed Hodge-Tate structure. We
calculate it and prove that:

\begin{theorem}\label{thm2}
The mixed Hodge structure given in Theorem \ref{main} for the
dilogarithm, that is,
$H^{2}(\mathcal{\overline{M}}_{0,5}-A(\vec{a}), B_2-B_2 \cap
A(\vec{a}))$, is isomorphic to the one given by P. Deligne. And
for our
case the period matrix is the following:\\
If $z\neq 0,1$, it equals
\[\left[
\begin{matrix}

 1         & 0               &0 \\
-Li_{1}(z) & 2\pi i          &0 \\
-Li_{2}(z) & 2\pi i\log{z}   &(2\pi i)^2

\end{matrix}
\right]\] If $z=1$, it is
\[\left[
\begin{matrix}

 1          & 0       \\
 -Li_{2}(1) &(2\pi i)^2

\end{matrix}
\right]\]
\end{theorem}

With above Theorem \ref{thm2}, the interesting question is that
for the classical $n$-th polylogarithm ($n\geq 3$), whether or not
our construction is isomorphic to the canonical one given by
Deligne. The general situation is more delicate than the
dilogarithm case. We have some results and it seems that they are
not isomorphic if $n\geq 3$. More detail will appear in
\cite{Wang2}.

\begin{remark}
Another construction of the multiple polylogarithm motives has
been given by A.Goncharov in \cite{Gon4} where he uses another
sequence of blowups. As \textbf{framed} mixed Tate motives, the
two constructions should be equivalent. But our construction is
canonical and more natural.
\end{remark}

\section{Geometry of the moduli space $\mathcal{\overline{M}}_{0,S}$ and $B_n$}
First recall that there is a one-to-one correspondence between the
boundary divisors of $\mathcal{\overline{M}}_{0,S}$ and the stable
unordered $2$-partitions of the set $S$. Let $\sigma=\sigma_1
|\sigma_{2}$ be a $2$-partition of the set $S$, then the stability
condition means that $|\sigma_1|\geq 2$ and $|\sigma_2|\geq 2$.
We'll denote by $D(\sigma)$ the corresponding boundary divisor.
\begin{definition}\label{def1}
Let $T=\{t_1<t_2< \dots, <t_k<t_1\}$ with the given cyclic order
$\rho$. A subset $A$ of $T$ is called strictly ordered if there
exists some $t_i\in T$ and $l$ a positive integer, such that
$A=\{t_{i}, t_{i+1},\dots, t_{i+l}\}$ (the subscripts are counted
mod $k$). That is, its elements are in consecutive order with
respect to $\rho$. Given a $2$-partition $\sigma$ of $T$,
$\sigma=\sigma_1 |\sigma_{2}$, we say that $\sigma$ is strictly
ordered with respect to $\rho$ if one of the $\sigma_i's$ is a
strictly ordered subset of $T$.
\end{definition}
For example, take $k=4$, then $A=\{s_1,s_2,s_3\}$ is strictly
ordered, but $B=\{s_1,s_3,s_4\}$ is not.

Now consider the open $n$-simplex
$\Delta_{n}=\{(t_{1},\dots,t_{n})\in\mathbb{R}^{n}|0<t_{1}<\dots<t_{n}<1\}$.
As mentioned in the Introduction, via cross-ratio,
$\mathcal{M}_{0,S}(\mathbb{C})$ is identified with the subset
$\{(x_1,\dots,x_n)\in \mathbb{C}^n|x_i\ne x_j,i\ne j\; ;x_k\ne
0,1; \,\  k=1,\dots,n\}$ of $\mathbb{C}^n$. Under this
identification, $\Delta_{n}$ is a subset of
$\mathcal{M}_{0,S}(\mathbb{C})$. Thus we have a natural map $\Phi$
which embeds $\Delta_{n}$ into
$\mathcal{\overline{M}}_{0,S}(\mathbb{C})$. Let $B_{n}$ be the
Zariski closure of the boundary of the closure of
$\Phi(\Delta_{n})$ in $\mathcal{\overline{M}}_{0,S}(\mathbb{C})$,
which is called algebraic Stasheff polytope (see \cite{GM}). Then
we have the following:

\begin{proposition}\label{prop1}
$B_{n}$ is a union of boundary divisors indexed by the stable
$2$-partitions of $S$ which are strictly ordered with respect to
the cyclic order $\rho$. That is, they correspond to breaking a
circle into two connected arcs. Furthermore, $B_{n}$ is an
``algebraic Stasheff polytope''. I.e., there is a bijection
between the irreducible components $D_{i}$ of $B_{n}$ and the
codimension one faces $F_{i}$ of the Stasheff polytope $K_{n+2}$
such that a subset of $D_{i}$'s has a non-empty intersection of
expected codimension if and only if the respective subset of
$F_{i}$'s has this property.
\end{proposition}

\begin{proof}
See \cite{GM}, Proposition 2.1.\qed
\end{proof}
\noindent For two unordered stable $2$-partitions $\sigma=\sigma_1
| \sigma_{2}$ and $\tau=\tau_1|\tau_2$ of S, we define:
\[a(\sigma,\tau):=\textit{the number of non-empty
intersections of}\  \sigma_i \cap \tau_j;\  i,j=1,2.\] Clearly,
$a(\sigma,\tau)=2,3$ or $4$, and $a(\sigma,\tau)=2$ if and only if
$\sigma=\tau$, which implies that the boundary divisors
$D(\sigma)=D(\tau)$.

\begin{lemma}\label{lemma1}
(1) $a(\sigma,\tau)=3$ if and only if $D(\sigma) \ne D(\tau)$ and
$D(\sigma) \cap D(\tau)\ne
\emptyset$.\\
(2) $a(\sigma,\tau)=4$ if and only if $D(\sigma) \cap D(\tau)=
\emptyset$.
\end{lemma}

\begin{proof}
(1) By definition, $a(\sigma,\tau)=3$ if and only if one of the
following holds: \[\sigma_1 \subsetneqq \tau_1,  \; \sigma_1
\subsetneqq \tau_2 , \;  \sigma_2 \subsetneqq \tau_1,  \; \sigma_2
\subsetneqq \tau_2.
\] By the Fact 4 in \cite[page~552]{Ke}, this is exactly the sufficient
and necessary condition for that $D(\sigma) \ne D(\tau)$ and
$D(\sigma) \cap D(\tau)\ne \emptyset$.

(2) By definition, $a(\sigma,\tau)=4$ means that there are
distinct elements $i, j, k, l \in S$ such that $i \in \sigma_1
\setminus \tau_1$, $j \in \tau_1 \setminus \sigma_1$, $k \in
\sigma_1 \cap \tau_1$, $l \notin \sigma_1 \cup \tau_1$. Then it
follows from the proof of the Fact 4 in \cite[page~552]{Ke}.
\qed
\end{proof}

\begin{definition}
A family of stable $2$-partitions $\{\sigma_1,\dots,\sigma_m\}$ of
$S$ is called good if $a(\sigma_i,\sigma_j)=3$, for $i\ne j$. A
family of boundary divisors $\{D(\sigma_1),\dots,D(\sigma_m)\}$ is
called compatible, if the corresponding partitions $\sigma_1,
\dots,\sigma_m$ are good.
\end{definition}

\begin{lemma}\label{lemma2}
Given $m$ pairwise distinct boundary divisors
$D(\sigma_1),\dots,D(\sigma_m)$, then
\[
 D(\sigma_1)\cap\dots\cap D(\sigma_m)\ne\emptyset
\]
if and only if $\{\sigma_1,\dots,\sigma_m\}$ is good.
\end{lemma}

\begin{proof}
First, if $\{\sigma_1,\dots,\sigma_n\}$ is not good, then
$a(\sigma_i,\sigma_j)\ne 3$ for some $i\ne j$. Since $D(\sigma_i)$
and $D(\sigma_j)$ are distinct, by Lemma \ref{lemma1},
$D(\sigma_i)\cap D(\sigma_j)=\emptyset$. Hence,
\[D(\sigma_1)\cap\dots\cap D(\sigma_m)=\emptyset.\]
Now suppose that $\{\sigma_1,\dots,\sigma_n\}$ is good. The by
\cite[Proposition~3.5.2]{Man2}, there exists a stable $S$-tree
$\tau$ with $m$ (internal) edges $\{e_1,\dots,e_m\}$ such that
$D(\sigma_i)=D(\sigma(e_i)), i=1,\dots,m$, where $\sigma(e_i)$ is
the stable $2$-partition of $S$ corresponding to the one edge
$S$-tree obtained by contracting all (internal) edges but $e_i$.
Denote by $D(\tau)$ the closure of the stratum parametrizing the
stable curves of the combinatorial type $\tau$. Then for each $i$,
we have $D(\tau)\subset D(\sigma(e_i))=D(\sigma_i)$. Hence,
$D(\sigma_1)\cap\dots\cap D(\sigma_m)\ne\emptyset$. And since the
boundary divisors meet transversely, we can conclude that in this
case $D(\sigma_1)\cap\dots\cap D(\sigma_m)=D(\tau)$.
\qed
\end{proof}
Now fix an $s_{i}\in S$ , and let $S^{\prime}=S-\{s_{i}\}$ with
the induced cyclic order
$\rho^{\prime}:0<s_{1}<\dots<s_{i-1}<s_{i+1}<\dots<s_{n}<1<\infty<0$.
Let
\[
\pi_{S^{\prime}}:\mathcal{\overline{M}}_{0,S}\longrightarrow
\mathcal{\overline{M}}_{0,S^{\prime}}
\]
be the contraction morphism which contracts the section marked by
$s_{i}$. Similarly, we have the natural embedding of
$\Delta_{n-1}=\{( t_{1},\dots,t_{i-1},t_{i+1},\dots,t_{n})
\in\mathbb{R}^{n-1}|0<t_{1}<\dots<t_{i-1}<t_{i+1}<\dots<t_{n}<1\}$
into $\mathcal{\overline{M}}_{0,S^{\prime}}(\mathbb{C})$ . And we
also have the corresponding algebraic Stasheff polytope $B_{n-1}$.
Now we have the following key property:

\begin{proposition}\label{prop3}
Under the contraction morphism $\pi_{S^{\prime}}$, the vertices of
$B_{n}$ are mapped to the vertices of $B_{n-1}$.
\end{proposition}
\begin{proof}
Let $v$ be a vertex of $B_n$. Then by Proposition \ref{prop1},
there exists exactly $n$ compatible boundary divisors of $B_n$,
say, $\{D(\sigma_1),\dots , D(\sigma_n)\}$ such that $v\in
D(\sigma_1)\cap \dots \cap D(\sigma_{n})$. Consider the images
$\pi_{S^{\prime}}(D(\sigma_j))$ in
$\mathcal{\overline{M}}_{0,S^{\prime}}(\mathbb{C})$, for
$j=1,2\dots, n$. There are two possible cases for the
$\sigma_j$'s:

Case 1: $\sigma_j=\{s_{i_0},
s_{i_0+1}\}|S\setminus\{s_{i_0},s_{i_0+1}\}$; or $\sigma_j=
\{s_{i_0-1}, s_{i_0}\}|S\setminus\{s_{i_0-1}, s_{i_0}\}$;
\\
Case 2: $\sigma_j=S_1|S_2$, with $|S_1|\geq 3$ and $s_{i_0}\in
S_1$.

In Case 1, we get
\[\pi_{S^{\prime}}(D(\sigma_j))=\mathcal{\overline{M}}_{0,S^{\prime}}.\]
and in Case 2, we have
\[\pi_{S^{\prime}}(D(\sigma_j))=D(\sigma_j^{\prime}),\] where
$\sigma_j^{\prime}=S_1\setminus \{s_{i_0}\}|S_2$ is a stable
partition of $S^{\prime}$, and $D(\sigma_j^{\prime})$ is an
irreducible component of $B_{n-1}$. It's clear that if
$\sigma_j\ne \sigma_k$ belong to Case 2 , then
$\pi_{S^{\prime}}(D(\sigma_j))\ne \pi_{S^{\prime}}(D(\sigma_k))$.
Now we claim that among the $\sigma_1$, $\dots$, $\sigma_m$, there
exists exactly one of them belonging to Case 1. Indeed, first
suppose all of them belong to Case 2. Then we have $n$ pairwise
distinct boundary divisors of
$\mathcal{\overline{M}}_{0,S^{\prime}}$ and their intersection is
empty due to the dimension consideration. But
$\pi_{S^{\prime}}(v)$ is contained in this intersection,
contradiction. Secondly, suppose there are two of them in Case 1.
Say they are $\sigma_1=\{s_{i_0},
s_{i_0+1}\}|S\setminus\{s_{i_0},s_{i_0+1}\}$ and
$\sigma_2=\{s_{i_0-1}, s_{i_0}\}|S\setminus\{s_{i_0-1},
s_{i_0}\}$. But then $a(\sigma_1, \sigma_2)=4$, thus
$D(\sigma_1)\cap D(\sigma_2)=\emptyset$. It contradicts that $v\in
D(\sigma_1)\cap \dots \cap D(\sigma_{n})$. Hence the claim is
proved.

By the claim, we see that $\pi_{S^{\prime}}(v)$ is contained in
the intersection of $n-1$ distinct compatible divisors of
$B_{n-1}$. By Proposition \ref{prop1}, $\pi_{S^{\prime}}(v)$ is a
vertex of $B_{n-1}$.
\qed
\end{proof}

Using Proposition \ref{prop3}, we can prove the following :

\begin{proposition}\label{prop4}
Let $\pi_{i}:\mathcal{\overline{M}}_{0,S}\longrightarrow
\mathcal{\overline{M}}_{0,\{0,s_{i},1,\infty\}}=\mathbb{P}^{1}$ be
the contraction morphism which forgets all sections but the ones
marked by $\{ 0,s_{i},1,\infty\}$. Then the images of the vertices
of $B_{n}$ are contained in the set $\{0,1\}$.
\end{proposition}
\begin{proof}
Under the identification of
$\mathcal{\overline{M}}_{0,\{0,s_{i},1,\infty\}}$ with
$\mathbb{P}^{1}$, $B_1$ is the interval $[0,1]$. Thus the vertices
are $0$, $1$. The map $\pi_{i}$ is a composition of $n-1$
contraction morphisms, each of which forgets only one section of
those labeled by the subset $S\setminus \{0,s_{i},1,\infty\}$. Now
we can apply Proposition \ref{prop3} repeatedly to those $n-1$
contraction morphisms. Therefore, the images of the vertices of
$B_{n}$ are either $0$ or $1$.
\qed
\end{proof}

Next, we'll introduce some non-boundary divisors in
$\mathcal{\overline{M}}_{0,S}$ and prove some properties which we
need later on. For each $i=1,\dots, n$, consider the morphism
\[\pi_{(n), i}:\mathcal{\overline{M}}_{0,S}\longrightarrow
\mathcal{\overline{M}}_{0,\{0,s_{i},1,\infty \}}=\mathbb{P}^{1}\]
which contracts all sections but those labeled by $0,
s_{i},1,\infty$. Now let $a\in \mathbb{P}^{1}\setminus
\{0,1,\infty \}$, that is, it lies in the open stratum
$\mathcal{M}_{0,\{0,s_{i},1,\infty\}}$. Then the inverse image
$\pi_{(n),i}^{-1}(a)$ is a divisor of
$\mathcal{\overline{M}}_{0,S}$. The following proposition
describes how $\pi_{(n),i}^{-1}(a)$ intersects the boundary
divisors of $\mathcal{\overline{M}}_{0,S}$.

\begin{proposition}\label{nonb1}
Let $\sigma=\sigma_1 | \sigma_{2}$ be a stable $2$-partition of
$S$ and $D(\sigma)$ be the corresponding boundary divisor.
Then:\\
(1) If $\{0,s_i \}$ and $\{1,\infty\}$, $\{1,s_i \}$ and
$\{0,\infty\}$, or $\{\infty,s_i \}$ and $\{1,0\}$ belong to
different parts of $\sigma$, then \[\pi_{(n),i}^{-1}(a)\cap
D(\sigma) =\emptyset.\] (2) Otherwise, we can assume that
$|\sigma_1\cap \{0,s_i,1,\infty \}|\geq 3$.
Then:\[\pi_{(n),i}^{-1}(a)\cap D(\sigma) =
\pi_{(|\sigma_1|-2),i}^{-1}(a)\times
\mathcal{\overline{M}}_{0,\sigma_2\cup \{u\}}\] where $u$ is some
fixed element of $\sigma_1$ and $\pi_{(|\sigma_1|-2),i}:
\mathcal{\overline{M}}_{0,\sigma_1\cup\{t\}}\longrightarrow
\mathcal{\overline{M}}_{0,\{0,s_{i},1,\infty \}}$ is the map
contracting all sections but those marked by $0, s_{i},1,\infty$;
and if $\{0,s_i,1,\infty \}\subset \sigma_1$, $t$ is some fixed
element of $\sigma_2$; otherwise, $|\sigma_1\cap \{0,s_i,1,\infty
\}|=3$, then $t$ is the only element of $\{0,s_i,1,\infty \}$ not
in $\sigma_1$.
\end{proposition}

\begin{proof}
For the first case, it's easy to check that the image of
$D(\sigma)$ under
$\pi_{(n),i}:\mathcal{\overline{M}}_{0,S}\longrightarrow
\mathcal{\overline{M}}_{0,\{0,s_{i},1,\infty \}}$ is contained in
the set $\{0,1,\infty\}$. But $a\neq 0,1,\infty$, thus
$\pi_{(n),i}^{-1}(a)\cap D(\sigma) =\emptyset.$

For the second case, consider the
morphism:\[\beta:\mathcal{\overline{M}}_{0,S}\longrightarrow
\mathcal{\overline{M}}_{0,\sigma_1\cup\{t\}}\times
\mathcal{\overline{M}}_{0,\sigma_2\cup \{u\}}\] which is the
product of contracting all sections labeled by $\sigma_2$ but $t$
and contraction of all sections labeled by $\sigma_1$ but $u$.

By the Fact 2 of \cite[page~551]{Ke}, the restriction of $\beta$
on $D({\sigma})$ is an isomorphism which is independent of the
choices of $u$ and $t$. Then we have the following commutative
diagram:
\[
\begin{CD}
D(\sigma)
@>\beta_{|D(\sigma)}>> \mathcal{\overline{M}}_{0,\sigma_1\cup\{t\}}\times\mathcal{\overline{M}}_{0,\sigma_2\cup \{u\}} \\
@V\pi_{(n),i|D(\sigma)}VV    @V pr_{1}VV \\
\mathcal{\overline{M}}_{0,\{0,s_{i},1,\infty \}}
@<\pi_{(|\sigma_1|-2),i}<<\mathcal{\overline{M}}_{0,\sigma_1\cup\{t\}}
\end{CD}
\]
where $\beta_{|D(\sigma)}$ and $\pi_{(n),i|D(\sigma)}$ are the
restrictions of the maps $\beta$ and $\pi_{(n),i}$ on the divisor
$D(\sigma)$ respectively, and $pr_1 $ is the projection on the
first factor.

Thus we obtain that:
\begin{align}
\notag
\pi_{(n),i}^{-1}(a)\cap D(\sigma) &= \pi_{(n),i|D(\sigma)}^{-1}(a)\\
\notag
                                 &= \beta_{|D(\sigma)}^{-1}(pr_1^{-1}(\pi_{(|\sigma_1|-2),i}^{-1}(a)))\\
\notag
                                 &= \beta_{|D(\sigma)}^{-1}(\pi_{(|\sigma_1|-2),i}^{-1}(a)\times\mathcal{\overline{M}}_{0,\sigma_2\cup \{u\}})
\notag
\end{align}
Since $\beta_{|D(\sigma)}$ is an isomorphism, we have
\[\pi_{(n),i}^{-1}(a)\cap D(\sigma)
=\pi_{(|\sigma_1|-2),i}^{-1}(a)\times
\mathcal{\overline{M}}_{0,\sigma_2\cup \{u\}}\] via the map
$\beta$. \qed
\end{proof}

Next we have the following description of $\pi_{(n),i}^{-1}(a)$:
\begin{proposition}\label{nonb2}
(1). $\pi_{(n),i}^{-1}(a)$ is an irreducible, reduced, smooth
divisor of $\mathcal{\overline{M}}_{0,S}$.\\
(2). It is a Tate variety, that's, its motive is a direct sum of
pure Tate motives.
\end{proposition}

\begin{proof}
We'll use the induction on $n$. When $n=1$, $S=\{0,s_1,1,\infty\}$
and
\[\pi_{(1),1}:\mathcal{\overline{M}}_{0,\{0,s_1,1,\infty\}}
\longrightarrow \mathcal{\overline{M}}_{0,\{0,s_1,1,\infty\}}\] is
the identity map. And the Proposition is obvious. Let's suppose
that it's true for $\pi_{(k),i}^{-1}(a)$, $k\leq n$, and consider
$\pi_{(n+1),i}^{-1}(a)$. Now $S=\{0,s_1,\dots,s_n, s_{n+1},
1,\infty\}$. Let $S^{\prime}=S\setminus \{s_{n+1}\}$. Then we have
the following commutative diagram:
\[
\begin{CD}
\mathcal{\overline{M}}_{0,S}@>\gamma>>\mathcal{\overline{M}}_{0,S^{\prime}}\times\mathcal{\overline{M}}_{0,\{0,s_{n+1},1,\infty\}}\\
@V\pi VV  @Vpr_1VV\\
\mathcal{\overline{M}}_{0,S^{\prime}}@=\mathcal{\overline{M}}_{0,S^{\prime}}\\
@V\pi_{(n),i}VV       @.\\
\mathcal{\overline{M}}_{0,\{0,s_i,1,\infty\}}
\end{CD}
\]
where $\pi$ is the map forgetting the section $s_{n+1}$, $\gamma$
is the product of $\pi$ and the morphism from
$\mathcal{\overline{M}}_{0,S}$ to
$\mathcal{\overline{M}}_{0,\{0,s_{n+1},1,\infty\}}$ contracting
all sections but those marked by $0,s_{n+1},1,\infty$; and $pr_1$
is the projection on the first factor. Moreover, we have:\[
\pi_{(n+1),i}=\pi_{(n),i}\circ \pi :
\mathcal{\overline{M}}_{0,S}\longrightarrow
\mathcal{\overline{M}}_{0,\{0,s_i,1,\infty\}}.\]
Hence,\[\pi_{(n+1),i}^{-1}(a)=\pi^{-1}(\pi^{-1}_{(n),i}(a)).\]
It's known that $\pi:\mathcal{\overline{M}}_{0,S}\longrightarrow
\mathcal{\overline{M}}_{0,S^{\prime}}$ is the universal curve on
$\mathcal{\overline{M}}_{0,S^{\prime}}$, so $\pi$ has
geometrically reduced fibers, and by
induction,$\pi^{-1}_{(n),i}(a)$ is reduced, thus
$\pi_{(n+1),i}^{-1}(a)$ is reduced.

In \cite{Ke}, Keel proved that the map $\gamma$ in the above
diagram is isomorphic to a sequence of blowups of
$\mathcal{\overline{M}}_{0,S^{\prime}}\times\mathcal{\overline{M}}_{0,\{0,s_{n+1},1,\infty\}}$
along all the boundary divisors of
$\mathcal{\overline{M}}_{0,S^{\prime}}$. From the above diagram,
we see that:
\begin{align}
\notag
\pi_{(n+1),i}^{-1}(a)&=\gamma^{-1}(pr_1^{-1}(\pi^{-1}_{(n),i}(a)))\\
\notag
                   &=\gamma^{-1}(\pi_{(n),i}^{-1}(a)\times\mathcal{\overline{M}}_{0,\{0,s_{n+1},1,\infty\}})
\notag \end{align} By the explicit blowups given in \cite{Ke},
$\pi_{(n+1),i}^{-1}(a)$ is the composition of blowups of
$\pi_{(n),i}^{-1}(a)\times\mathcal{\overline{M}}_{0,\{0,s_{n+1},1,\infty\}}$
along the intersections of $\pi_{(n),i}^{-1}(a)$ with the boundary
divisors. By the Proposition \ref{nonb1}, these intersections are
either empty or of the form: \[\pi_{(l),i}^{-1}(a)\times
\mathcal{\overline{M}}_{0,T}\] where $l<n$ and $T\subset S$.

By the induction assumption, $\pi_{(l),i}^{-1}(a)\times
\mathcal{\overline{M}}_{0,T}$ is irreducible and smooth, and
$\pi_{(n),i}^{-1}(a)\times\mathcal{\overline{M}}_{0,\{0,s_{n+1},1,\infty\}}$
is an irreducible smooth divisor of
$\mathcal{\overline{M}}_{0,S^{\prime}}\times\mathcal{\overline{M}}_{0,\{0,s_{n+1},1,\infty\}}$.
Altogether we see that $\pi_{(n+1),i}^{-1}(a)$ is the composition
of blowups of an irreducible smooth variety along an irreducible
smooth subvariety. Therefore, it is irreducible and smooth. And
this proves the first part of the proposition.

For the second part, we need the following formula of the motive
of the blowup of a variety along a subvariety. This was proved by
Manin in \cite[Corollary, page~463]{Man1}. Let $X^{\prime}$ be the
blowup of a variety $X$ along a subvariety $Y$ of codimension $r$.
Then we have:
\[h(X^{\prime})=h(X)\oplus(\oplus_{i=1}^{r-1}h(Y)\otimes\mathbb{L}^i).\]
where $h(X^{\prime})$, $h(X)$ and $h(Y)$ denote the motives of
$X^{\prime}$, $X$ and $Y$ respectively; $\mathbb{L}$ is the Tate
motive, $h(\mathbb{P}^1)$.

Therefore, if both $X$ and $Y$ are Tate varieties, by the formula
above, the blowup $X^{\prime}$ is also a Tate variety. Because the
product of Tate varieties is also a Tate variety, and by
induction, $\pi_{(n+1),i}^{-1}(a)$ is the composition of blowups
of a Tate variety along a Tate subvariety. Thus it's a Tate
variety. This proves the second part.
\qed
\end{proof}

Now we consider the intersections of the divisors
$\pi_{(n),i}^{-1}(a)$. Let $a_1, a_2, \dots a_m\in
\mathbb{C}\setminus\{0,1\}$, and $m\leq n$. Then we have:
\begin{proposition}\label{nonb3}
The intersection $\pi_{(n),1}^{-1}(a_1)\cap
\pi_{(n),2}^{-1}(a_2)\cap \cdots \cap\pi_{(n),m}^{-1}(a_m)$ is a
smooth irreducible subvariety; it's also a Tate variety.
\end{proposition}
\begin{proof}
First let's look at the case $m=2$. Notice for any irreducible
boundary divisor $D(\sigma)$, where $\sigma$ is some stable
$2$-partition of $S$, we have \[\pi_{(n),1}^{-1}(a_1)\cap
\pi_{(n),2}^{-1}(a_2)\cap
 D(\sigma)=(\pi_{(n),1}^{-1}(a_1)\cap D(\sigma))\cap(\pi_{(n),2}^{-1}(a_2)\cap
 D(\sigma)).\] By the Proposition \ref{nonb1},
$\pi_{(n),1}^{-1}(a_1)\cap \pi_{(n),2}^{-1}(a_2)\cap D(\sigma)$ is
either empty or of the form:
\[(\pi_{(|\sigma_1|-2),1}^{-1}(a_1)\cap
\pi_{(|\sigma_1|-2),2}^{-1}(a_2)) \;\times
\mathcal{\overline{M}}_{0,\sigma_2\cup \{u\}}.\] Here we use the
notations in Proposition \ref{nonb1}. Now by induction on $n$, we
see that $\pi_{(n),1}^{-1}(a_1)\cap \pi_{(n),2}^{-1}(a_2)$ is
obtained by the blowups of smooth irreducible Tate varieties along
the Tate subvarieties. Hence it's also smooth, irreducible and
Tate.

When $m>2$, similarly for any irreducible boundary divisor
$D(\sigma)$, we have
\[(\bigcap_{i=1}^{m}\pi_{(n),i}^{-1}(a_i))\bigcap
D(\sigma)=\bigcap_{i=1}^{m}(\pi_{(n),i}^{-1}(a_i)\cap
D(\sigma)).\] Then similarly, use Proposition \ref{nonb1} and
induction on $n$, it's true for any $m$. \qed
\end{proof}

\begin{corollary}\label{nonb4}
Let $\sigma_1, \sigma_2, \dots, \sigma_k$ be stable $2$-partitions
of $S$. For each $i=1,\dots, k$, $D(\sigma_i)$ is the
corresponding boundary divisor. Then the intersection
\[\pi_{(n),1}^{-1}(a_1)\cap \pi_{(n),2}^{-1}(a_2)\cap \cdots
\cap\pi_{(n),m}^{-1}(a_m)\cap D(\sigma_1)\cap D(\sigma_2)\cap
\dots \cap D(\sigma_k)\] is either empty or a smooth irreducible
Tate variety.
\end{corollary}
\begin{proof}
By the combinatorial description of the intersection
$D(\sigma_1)\cap D(\sigma_2)\cap \dots \cap D(\sigma_k)$, (for
$k=2$, see \cite[Fact~4]{Ke}, and for general $k$, it can be done
inductively) we see that the non-empty intersection in the
corollary can be written as the product of irreducible smooth Tate
varieties. Thus, it's also an irreducible smooth Tate variety.
\qed
\end{proof}

\section{The meromorphic form $\Omega_{S}(\vec{a})$ and The divisor $A_S(\vec{a})$}
\subsection{The form $\Omega_{S}(\vec{a})$}\,

Let $\vec{a}=(a_{s_1},\cdots,a_{s_n})\in \mathbb{C}^n$ and
consider the following $n$ contraction morphisms:
\[\beta_{s_i}:\mathcal{\overline{M}}_{0,S}\longrightarrow
\mathcal{\overline{M}}_{0,\{0,s_i,1,\infty\}}=\mathbb{P}^1 \qquad
i=1,\dots, n\] which forget all sections but those labeled by
$0,s_i,1,\infty$. We define the meromorphic $n$-form
$\Omega_{S}(\vec{a})$ of $\mathcal{\overline{M}}_{0,S}$ as
follows:
\begin{equation}
\Omega_{S}(\vec{a}):=\frac{d\beta_{s_1}}{\beta_{s_1}-a_{s_1}}\wedge
...\wedge\frac{d\beta_{s_n}}{\beta_{s_n}-a_{s_n}}\label{form1}\\
\end{equation}
Furthermore, we have another useful description of
$\Omega_{S}(\vec{a})$. Let
\[\beta:\mathcal{\overline{M}}_{0,S}\longrightarrow
\prod_{i=1}^{n}\mathcal{\overline{M}}_{0,\{0,s_i,1,\infty\}}=(\mathbb{P}^1)^n\]
be the birational proper morphism of the products of
$\beta_{s_i}$, $i=1,\dots,n$. And let $(t_{s_1},\dots,t_{s_n})$ be
the affine coordinates of $(\mathbb{P}^1)^n$, then by definition,
we have
\[\Omega_{S}(\vec{a})=\beta^{*}(\frac{dt_{s_1}}{t_{s_1}-a_{s_1}}\wedge
...\wedge\frac{dt_{s_n}}{t_{s_n}-a_{s_n}}).\]
\subsection{The divisor $A_S(\vec{a})$}
We define $A_{S}(\vec{a})$ to be the divisor of singularities of
$\Omega_{S}(\vec{a})$. Following \cite{GM}, for $\alpha
\in\{0,1,\infty\}$, we define $S(\alpha)$ by
\[ S(0):=\{s\in S|a_{s}=0\}, \ S(1):=\{s\in S |a_{s}=1\}\]
and
\[S(0,1):=S(0)\cup S(1), \ S(\infty):=S\setminus \{0,1,\infty\}.
\]
we say that a $2$-partition of $S$ has \textit{type $\alpha$} if
one part of it is of form $\{\alpha\}\cup T$ where $T$ is a
non-empty subset of $S(\alpha)$. Now we have the explicit
description of the divisor $A_{S}(\vec{a})$:

\begin{proposition}\label{divisor-a}
The divisor $A_{S}(\vec{a})$ of singularities of
$\Omega_S(\vec{a})$ on $\mathcal{\overline{M}}_{0,S}$ consists of
the following
irreducible components: \\
(a) The boundary divisors $D(\sigma)$ corresponding to those
stable $2$-partitions of $S$ which have some type $\alpha$,
$\alpha \in\{0,1,\infty\}$;\\
(b) The non-boundary divisors $\pi_{(n),i}^{-1}(a_{s_i})$, for
each $s_i\notin S(0,1)$; where $\pi_{(n),i}:
\mathcal{\overline{M}}_{0,S}\longrightarrow
\mathcal{\overline{M}}_{0,\{0,s_i,1,\infty\}}$ is the morphism
contracting all sections but those labeled by  $0,s_i,1,\infty$.
Moreover, $A_{S}(\vec{a})$ is a normal-crossing divisor.
\end{proposition}
\begin{proof}
We'll prove it by induction on $n$. For the case $n=1$,
$S=\{0,s_1,1, \infty\}$ and
$\mathcal{\overline{M}}_{0,S}=\mathbb{P}^1$. Then
$\Omega_{S}(\vec{a})=\frac{dt}{t-a_{s_1}}$, where $t$ is the
affine coordinate of $\mathbb{P}^1$. Thus
$A_{S}(\vec{a})=(a_{s_i})+(\infty)$, and the proposition is clear
in this case.

Assume that it's true for $n$. Now $S=\{0,s_1,\dots, s_n, s_{n+1},
1,\infty\}$ and $\vec{a}=(a_{s_1},\dots, a_{s_n}, a_{s_{n+1}})$.
Let $\vec{a}^{\prime}=(a_{s_1},\dots, a_{s_n})$,
$S^{\prime}=S\setminus \{s_n\}$, and
$A_{S^{\prime}}(\vec{a}^{\prime})$ the divisor of singularities of
the meromorphic form $\Omega_{S^{\prime}}(\vec{a}^{\prime})$ in
$\mathcal{\overline{M}}_{0,S^{\prime}}$. Consider the morphism
$\beta$ of the product of two contraction maps:
\[\beta:\mathcal{\overline{M}}_{0,S}\longrightarrow
\mathcal{\overline{M}}_{0,S^{\prime}}\times
\mathcal{\overline{M}}_{0,\{ 0,s_{n+1},1,\infty \}}.\] Then we
obtain that
\[\Omega_{S}(\vec{a})=\beta^{*}(\Omega_{S^{\prime}}(\vec{a}^{\prime})\wedge
\frac{dt}{t-a_{s_{n+1}}})\] Let $A(a_{s_{n+1}})$ be the divisor of
singularities of the form $\frac{dt}{t-a_{s_{n+1}}}$ in
$\mathcal{\overline{M}}_{0,\{ 0,s_{n+1},1,\infty \}}$. Clearly,
the divisor of singularities of the meromorphic form
$\Omega_{S^{\prime}}(\vec{a}^{\prime})\wedge
\frac{dt}{t-a_{s_{n+1}}}$ in
$\mathcal{\overline{M}}_{0,S^{\prime}}\times
\mathcal{\overline{M}}_{0,\{ 0,s_{n+1},1,\infty \}}$ is
$pr_1^{*}(A_{S^{\prime}}(\vec{a}^{\prime}))+pr_2^{*}(A(a_{s_{n+1}}))$,
where $pr_1$ and $pr_2$ are the projections on the first and
second factors respectively. And it's a normal crossing divisor.
By \cite{Ke}, $\beta$ is isomorphic to a sequence of blowups of
$\mathcal{\overline{M}}_{0,S^{\prime}}\times\mathcal{\overline{M}}_{0,\{0,s_{n+1},1,\infty\}}$
along all the boundary divisors of
$\mathcal{\overline{M}}_{0,S^{\prime}}$. By the following Lemma
\ref{blowup} on blowups, and Proposition \ref{nonb1} and
\ref{nonb2}, we can conclude that
$\beta^{*}(pr_1^{*}(A_{S^{\prime}}(\vec{a}^{\prime}))+pr_2^{*}(A(a_{s_{n+1}})))$
is a normal crossing divisor in $\mathcal{\overline{M}}_{0,S}$.
$A_{S}(\vec{a})$ is contained in it. Thus it's also a normal
crossing divisor.

On the other hand, the meromorphic form $\Omega_{S}(\vec{a})$ does
not necessarily have as singularities all the exceptional divisors
of $\beta$ in
$\beta^{*}(pr_1^{*}(A_{S^{\prime}}(\vec{a}^{\prime}))+pr_2^{*}(A(a_{s_{n+1}})))$
. In fact, $A_{S}(\vec{a})$ equals
$\beta^{*}(pr_1^{*}(A_{S^{\prime}}(\vec{a}^{\prime}))+pr_2^{*}(A(a_{s_{n+1}})))$
minus those spurious divisors.

By \cite[Lemma~3.8]{Gon4}, the spurious divisors are those
irreducible boundary divisors in $\mathcal{\overline{M}}_{0,S}$
which get blown down by $\beta$ to a subvariety of the product
which is not a stratum of the divisor
$pr_1^{*}(A_{S^{\prime}}(\vec{a}^{\prime}))+pr_2^{*}(A(a_{s_{n+1}}))$.

By \cite[Lemma~1, Page~554]{Ke}, the exceptional divisors of
$\beta$ are exactly those boundary divisors corresponding to the
following $2$-stable partitions:

$(*):\quad \sigma=\sigma_1|\sigma_2$, with $s_{n+1}\in \sigma_1$,
$|\sigma_1|\geq 3$, and $|\sigma_1\cap \{0,1,\infty\}|\leq 1$.

By induction, the divisor $A_{S^{\prime}}(\vec{a}^{\prime})$
consists of the following four types of irreducible divisors (we
identify the boundary divisors with the corresponding $2$-stable
partitions here):
\[  0T_0|\cdots1\infty ; \qquad  1T_1|\cdots 0\infty ; \qquad  \infty T_{\infty}|\cdots
01 ;\qquad  \pi_{(n),i}^{-1}(a_{s_i}), \,\,\; a_{s_i}\ne 0,1\]
where $T_{\alpha}$ is a non-empty subset of $S(\alpha)$, $\alpha
\in \{0,1,\infty\}$.

First, notice the pullback of the divisor
$\pi_{(n),i}^{-1}(a_{s_i})$ equals $\pi_{(n+1),i}^{-1}(a_{s_i})$
which is not an exceptional divisor of $\beta$, thus it's a
component of the divisor $A_{S}(\vec{a})$.

For each $\alpha\in\{0,1,\infty\}$, by \cite[Fact~3,
page~552]{Ke}, the pullback of the boundary divisor $\alpha
T_{\alpha}|\cdots$ in the list above is equal to the sum of the
following two components:
\[\alpha T_{\alpha}s_{n+1}|\cdots\;+\; \alpha T_{\alpha}|s_{n+1}\cdots.\]

By the condition $(*)$, $\alpha T_{\alpha}|s_{n+1}\cdots$ is not
an exceptional divisor of $\beta$, hence it's a component of the
divisor $A_{S}(\vec{a})$; but $\alpha T_{\alpha}s_{n+1}|\cdots$ is
an exceptional divisor. We need to check if it's a spurious
divisor. Notice that \[\beta(\alpha
T_{\alpha}s_{n+1}|\cdots)=\alpha T_{\alpha}|\cdots \; \times \;
\alpha s_{n+1}|\{0,1,\infty\}\setminus \alpha\] Now clearly
$\alpha T_{\alpha}|\cdots\; \times\; \alpha
s_{n+1}|\{0,1,\infty\}\setminus \alpha$ is a stratum of the
divisor
$pr_1^{*}(A_{S^{\prime}}(\vec{a}^{\prime}))+pr_2^{*}(A(a_{s_{n+1}}))$
if and only if $\alpha$ is in the divisor $A(a_{s_{n+1}})$,
therefore, $\alpha T_{\alpha}s_{n+1}|\cdots$ is a component of the
divisor $A_{S}(\vec{a})$ if and only if $\alpha$ is in the divisor
$A(a_{s_{n+1}})$. For $\alpha=0,1$, it depends upon whether or not
$a_{s_{n+1}}=\alpha$. And $\alpha=\infty$ is always in the divisor
$A(a_{s_{n+1}})$, hence $\infty T_{\infty}s_{n+1}|\cdots 01$ is a
component of $A_{S}(\vec{a})$.

Now there are only two more components of $A_{S}(\vec{a})$. One is
$\infty s_{n+1}|\cdots 01$. It's not an exceptional divisor and is
contained in the pullback of $\infty s_{n+1}|01$. The other one is
$\pi_{(n+1),n+1}^{-1}(a_{s_{n+1}})$, when $a_{s_{n+1}}\ne 0,1$.
And it is contained in the pullback of $A(a_{s_{n+1}})$. Thus we
have all the components of $A_{S}(\vec{a})$ listed in the
proposition. This concludes the proof.
\qed
\end{proof}
Now let's prove a lemma on blowups which is used in the proof of
the preceding proposition.
\begin{lemma}\label{blowup}
Let $M$ be a smooth complex variety of dimension $n$, $B$ an
irreducible smooth subvariety of $M$ with codimension $m \geq 2$,
and $D=D_1\cup D_2\cup \dots \cup D_k$ a normal-crossing divisor
of $M$ with smooth irreducible components $D_i$, $i=1,\dots k$.
Suppose that for each $D_i$, only one of the following cases
happens: $B\subset D_i$, $B\cap D_i=\emptyset$, or: $B\cap D_i$ is
an irreducible smooth subvariety with codimension $m-1$. Let
$\pi:M_{B}^{\prime}\longrightarrow M$ be the blowup of $M$ along
$B$. Let $D^{\prime}_{i}$ be the strict transform of $D_i$, $1\leq
i\leq k$, and $E=\pi^{-1}(B)$ the exceptional divisor. Then the
divisor $D^{\prime}=D^{\prime}_1\cup D^{\prime}_2\cup \dots \cup
D^{\prime}_k\cup E$ is a normal-crossing divisor in
$M_{B}^{\prime}$.
\end{lemma}

\begin{proof}
Let $x^{\prime}\in D^{\prime}$. If $x^{\prime}\notin E$, the $\pi$
is locally at $x^{\prime}$ an isomorphism. So we only need to
consider the case $x^{\prime}\in E$. Then $x=\pi(x^{\prime})\in
B$. By rearranging the indices of the $D_i$, if necessary, we may
assume that $B\subset D_i$ for $1\leq i\leq a$, $x\in D_j\cap
B\neq \emptyset$ for $a+1\leq j\leq b$, and $x\notin D_j\cap B$
for $j>b$. Here $1\leq a\leq \min\{m,k\}$, $a\leq b\leq
\min\{k,m\}$ and $b-a\leq n$. Since $D$ is a normal-crossing
divisor, there exists an open neighborhood $U$ of $x$ with local
coordinates $(z_1,z_2,\dots, z_n)$ such that
$z_1(x)=z_2(x)=\dots=z_n(x)=0$ and $U\cap B=\{y\in
U|z_1(y)=z_2(y)=\dots =z_m(y)=0\}$, for $1\leq i\leq a$, $U\cap
D_i=\{y\in U|z_i(y)=0\}$; for $a+1\leq j\leq b$, $U\cap D_j=\{y\in
U|z_j(y)=0\}$. And $\pi^{-1}(U)=\{(y,l)\in U\times
\mathbb{P}^{m-1}|z_i(y)l_j=z_j(y)l_i, 1\leq i,j\leq m\}$. Here
$l=[l_1:\dots:l_m]$ is the homogenous coordinate of
$\mathbb{P}^{m-1}$. Then a point $(y,l)$ in the inverse image of
$\pi^{-1}(U\cap(D_i\setminus B))$ has the property that $l_i=0$.
Now in $D^{\prime}$, suppose that $x^{\prime}\in D^{\prime}_i$ for
$1\leq i\leq c$, and $x^{\prime}\notin D^{\prime}_j$ for $j>c$.
Then $c\leq b$. Now since $x^{\prime}=(y,l)\in \pi^{-1}(U)$, it
can not happen that all $l_i=0$. Then there exists some $i_0>c$
such that $l_{i_0}\ne 0$. Consider the open neighborhood
$U_{i_0}=\{(y,l)\in\pi^{-1}(U)|l_{i_0}\ne 0\}$ of $\pi^{-1}(U)$.
And the restriction of the blowup $\pi$ on the $U_{i_0}$ has the
following
formula:\[(s_1,\cdots,z_{i_0},\cdots,s_m,z_{m+1},\cdots,z_n)
\mapsto (s_1z_{i_0},\cdots,z_{i_0},\cdots,s_mz_{i_0},
z_{m+1},\cdots,z_n)\] where $s_j=l_j/l_{i_0}$, $1\leq j\leq
m$,$j\ne i_0$. Therefore for $1\leq j\leq c$, $D^{\prime}_j\cap
U_{i_0}=\{(y,l)\in U_{i_0}|s_j=0\}$ and $E\cap U_{i_0}=\{(y,l)\in
U_{i_0}|z_{i_0}(y)=0\}$. Thus $D^{\prime}$ is a normal-crossing
divisor in $M_{B}^{\prime}$. \qed
\end{proof}

Next we'll prove the main theorem in this section:
\begin{theorem}\label{thm1}
The divisor $A_S(\vec{a})$ of singularities $\Omega_S(\vec{a})$
does not contain any $k$-dimensional face of the algebraic
stasheff polytope $B_{n}$, $0\leq k\leq n$ . Here
$\vec{a}=(a_{s_1},\dots,a_{s_n})$ with $a_{s_1}\ne 0$ and
$a_{s_n}\ne 1$.
\end{theorem}
\begin{proof}{Proof.}
First observe that if $A_S(\vec{a})$ contains some $k$-dimensional
face of $B_{n}$, then it also contains the vertices of this face.
Therefore it suffices to show that $A_S(\vec{a})$ does not contain
any vertex of $B_{n}$.

We'll prove that each irreducible component of $A_S(\vec{a})$
listed in the preceding Proposition \ref{divisor-a} can't contain
any vertex of $B_{n}$.

First consider the non-boundary divisor
$\pi_{(n),i}^{-1}(a_{s_i})$, for $a_{s_i}\ne 0,1$. Clearly, its
image under the contraction morphism: \[\pi_{(n),i}:
\mathcal{\overline{M}}_{0,S}\longrightarrow
\mathcal{\overline{M}}_{0,\{0,s_i,1,\infty\}}\] is $a_{s_i}$. On
the other hand, by the Proposition \ref{prop4} in Section 2, the
image of a vertex of $B_n$ under $\pi_{(n),i}$ is $0$, or $1$.
Since $a_{s_i}\ne 0,1$, $\pi_{(n),i}^{-1}(a_{s_i})$ doesn't
contain any vertex of $B_n$.

Let's consider the boundary divisor in $A_S(\vec{a})$. It has one
of the following types \[  0T_0|\cdots1\infty ; \qquad 1T_1|\cdots
0\infty ; \qquad  \infty T_{\infty}|\cdots 01 \] By the
Proposition \ref{prop1} in Section 2, we know that a vertex of
$B_n$ is contained in the intersection of exactly $n$ irreducible
components of $B_n$. Notice that the intersection of any $n+1$
irreducible boundary divisors of $\mathcal{\overline{M}}_{0,S}$ is
empty, therefore, it's enough to show that none of the boundary
divisors in $A_S(\vec{a})$ appears as an irreducible component of
$B_n$.

By the Proposition \ref{prop1} in Section 2, the irreducible
components of $B_n$ correspond to the stable $2$-partitions of $S$
which are strictly ordered with respect to the cyclic order $\rho:
0<s_1<s_2<\dots<s_n<1<\infty<0$. The partition
$0T_0|\cdots1\infty$ is not strict with respect to $\rho$ because
$a_{s_1}\ne 0$ then $s_1$ and $\infty$ separate any element of
$T_0$ from $0$; similarly, $1T_1|\cdots 0\infty$ is not strict
with respect to $\rho$ because $a_{s_n}\ne 1$ and $s_n$ and
$\infty$ block any element of $T_1$ to $1$; finally $\infty
T_{\infty}|\cdots 01$ is not strict with respect to $\rho$ because
$1$ and $0$ separate $\infty$ from any element of $T_{\infty}$.
Hence, none of the boundary divisors in $A_S(\vec{a})$ contains a
vertex of $B_n$. This concludes the proof.
\qed
\end{proof}

\section{Multiple polylogarithm motives}
We'll continue to use the notations in the previous sections. Now
consider the convergent iterated integral:
\begin{equation}\label{ite1}
I_{\gamma}(a_{s_1},\cdots,a_{s_n})=\int_{\gamma(\Delta_{n})}\frac{dt_{s_1}}{t_{s_1}-a_{s_1}}\wedge
\cdots \wedge\frac{dt_{s_n}}{t_{s_n}-a_{s_n}} \quad
a_{s_1}\neq0,a_{s_n} \neq1.\\
\end{equation}
where $\gamma:[0,1] \rightarrow \mathbb{C}$ is a piecewise smooth
simple path from $0$ to $1$ and $a_{s_i}\notin \gamma((0,1))$,
$i=1,\dots,n$,
$\Delta_{n}=\{(t_{s_1},\cdots,t_{s_n})\in\mathbb{R}^{n}|
0<t_{s_1}<\cdots<t_{s_n}<1\}$ is an open $n$-simplex in
$\mathbb{R}^{n}$, and
$\gamma(\Delta_{n})=\{(\gamma(t_{s_1}),\cdots,\gamma(t_{s_n}))|(t_{s_1},\cdots,t_{s_n})\in
\Delta_{n}\}$.

As mentioned in the Introduction, we can identify the open stratum
$\mathcal{M}_{0,S}(\mathbb{C})$ with the set $\{(x_1,\dots,x_n)\in
\mathbb{C}^n|x_i\ne x_j,i\ne j\; ;x_k\ne 0,1 \,\, ,
k=1,\dots,n\}$. Then $\gamma(\Delta_{n})$ is a subset of
$\mathcal{M}_{0,S}(\mathbb{C})$, thus we have a natural map
$\Phi_n$ which embeds $\gamma(\Delta_{n})$ into
$\mathcal{\overline{M}}_{0,S}(\mathbb{C})$.

Clearly we have the following equality:
\begin{equation}\label{pe}
I_{\gamma}(a_{s_1},\cdots,a_{s_n})=\int_{\Phi(\gamma(\Delta_{n}))}\Omega_{S}(\vec{a})
\end{equation}
where $\Omega_{S}(\vec{a})$ is the meromorphic $n$-form in
$\mathcal{\overline{M}}_{0,S}$ defined in Section 2.

Let $\Phi_n(\gamma)$ be the closure of $\gamma(\Delta_{n})$ in
$\mathcal{\overline{M}}_{0,S}$. Then the Zariski closure of the
boundary of $\Phi_n(\gamma)$ is $B_n$. This follows from the fact
that $B_n$ is the Zariski closure of the boundary of the closure
of $\Delta_n$ in $\mathcal{\overline{M}}_{0,S}$. Hence we have a
relative homology class:\[[\Phi_n(\gamma)]\in
H_n(\mathcal{\overline{M}}_{0,S}(\mathbb{C}),B_n(\mathbb{C});\mathbb{Q})\]
The meromorphic $n$-form $\Omega_{S}(\vec{a})$ gives rise to a
cohomology class: \[[\Omega_{S}(\vec{a})]\in
H^{n}_{DR}(\mathcal{\overline{M}}_{0,S}(\mathbb{C})-A_S(\vec{a})(\mathbb{C});
\mathbb{C}).\]

Now given the iterated integral (\ref{ite1}), we'll define:
\[
I^{\mathcal{M}}(a_{s_1},...,a_{s_n}):=H^{n}(\mathcal{\overline{M}}_{0,S}-A_S({\vec{a}}),
B_n-B_n\cap A_S({\vec{a}})).\]  Then we have the following main
theorem:

\begin{theorem}\label{mt}
$I^{\mathcal{M}}(a_{s_1},...,a_{s_n})$ carries an $n$-framed
Hodge-Tate structure with the frames coming from
$\Omega_S{(\vec{a})}$ and $\Phi_n(\gamma)$. The period is just the
iterated integral $I_{\gamma}(a_{s_1},...,a_{s_n})$. Moreover, if
all the $a_{s_i}$ are elements of a number field $F$, then
$I^{\mathcal{M}}(a_{s_1},..., a_{s_n})$ is a framed mixed Tate
motive over $F$.
\end{theorem}

Before we prove the theorem, first let's briefly recall the
definitions of framed Hodge-Tate structure and its period.
\begin{definition}
A mixed Hodge $\mathbb{Q}$-structure consists of the following
data:\\
1) a finite-dimensional $\mathbb{Q}$-vector space $H_{\mathbb{Q}}$
with a finite
increasing filtration $W_n$ called the weight filtration;\\
2) a finite decreasing filtration $F^p$ of
$H_{\mathbb{C}}:=H_{\mathbb{Q}}\otimes_{\mathbb{Q}} \mathbb{C}$
called the Hodge filtration;\\ these data satisfy the following
conditions:\\
For each associated graded piece
$Gr_{k}^{W}(H_{\mathbb{C}})=\frac{W_nH_{\mathbb{Q}}\otimes
\mathbb{C}}{W_{n-1}H_{\mathbb{Q}}\otimes \mathbb{C}}$, one has the
decomposition:\[Gr_{k}^{W}(H_{\mathbb{C}})=\bigoplus_{p+q=k}H^{p,q}\]
with \[H^{p,q}=F^pGr_{k}^{W}(H_{\mathbb{C}})\cap
\overline{F}^qGr_{k}^{W}(H_{\mathbb{C}}), \; \text{and}\;
H^{p,q}=\overline{H^{q,p}}.\] Here "\textemdash"  means the
complex conjugation of $H_{\mathbb{C}}$ with respect to
$H_{\mathbb{R}}:= H_{\mathbb{Q}}\otimes_{\mathbb{Q}} \mathbb{R}$.
The Hodge numbers are the integers
$h^{p,q}=dim_{\mathbb{C}}H^{p,q}=h^{q,p}$.
\end{definition}
By definition, a \emph{Hodge-Tate structure} $H$ is a mixed Hodge
$\mathbb{Q}$-structure $H$ with the Hodge numbers $h^{p,q}=0$
unless $p=q$. This means that for the weight filtration,
$Gr_{2n+1}^{W}H=0$ and $Gr_{2n}^{W}H$ is a finite direct sum of
$\mathbb{Q}(-n)$. It also implies that, for each $p\in
\mathbb{Z}$, the natural map
\begin{equation}\label{ht}
F^pH_{\mathbb{C}}\cap W_{2p}H_{\mathbb{C}}\rightarrow
Gr_{2p}^{W}H_{\mathbb{C}}
\end{equation} is an isomorphism.
\begin{definition}
An \emph{$n$-framed Hodge-Tate structure} $H$ is a Hodge-Tate
structure $H$ equipped with a nonzero vector in $Gr^{W}_{2n}H$ and
a nonzero functional on $Gr^{W}_{0}H$, that is, we have two
nonzero morphisms: \[v: \mathbb{Q}(-n)\rightarrow Gr^{W}_{2n}H, \
f: \mathbb{Q}(0)\rightarrow Gr^{W}_{0}(H)^{\vee}.\]
\end{definition}
To define the period of an $n$-framed Hodge-Tate structure, we
need to choose a map of $\mathbb{Q}$-vector spaces $F :
\mathbb{Q}\rightarrow H_{\mathbb{Q}}^{\vee}$ which lifts $f$, that
is, $Gr^{W}_{0}F=f$. Now let $f^{\prime}=F(1)$. Consider the
composition: \[\mathbb{Q}(-n)\rightarrow
Gr^{W}_{2n}H_{\mathbb{Q}}\rightarrow F^nH_{\mathbb{C}}\cap
W_{2n}H_{\mathbb{C}},\] where the first one is $v$ and the second
is provided by (\ref{ht}). It gives rise to a vector
$v^{\prime}\in F^nH_{\mathbb{C}}\cap W_{2n}H_{\mathbb{C}}$. The
period is the number $\langle v^{\prime}, f^{\prime}\rangle$. A
different choice of the lifting $F$ will change this period by
$2\pi i \times$ "weight $n-1$ period". In this sense, periods are
multi-valued.

For more information about periods of framed mixed Tate Motives or
periods of framed Hodge-Tate structures, we refer to
\cite[Chapter~5]{Gon2} and \cite[Section~3.2]{Gon3}.

Now let's prove some lemmas about Hodge-Tate structures.

\begin{lemma}\label{h-t1}
The category of Hodge-Tate structures is abelian. It's closed
under subquotients and extensions.
\end{lemma}
\begin{proof}
It's straightforward (c.f. \cite[Th\'{e}or\`{e}m~2.3.5]{De2}).
\qed
\end{proof}

\begin{lemma}\label{h-t3}
Let $A$ be a complex algebraic variety. Suppose that
\[A=\bigcup_{i=1}^{k}{A_i},\] where $A_i$ are closed subvarieties of $A$,
$i=1,\dots,k$. Let $d$ be a fixed positive integer. For any
non-empty subset $I\subset \{1,2,\cdots, k\}$, denote by
$A_{I}=\cap_{i \in I}A_i$. Suppose that all the cohomology groups
of $A_I$ carry the Hodge-Tate structures for all $I$.  Then
$H^d(A, \mathbb{Q})$ also carries a Hodge-Tate structure.
\end{lemma}
\begin{proof}
We'll use the Mayer-Vietoris exact sequence and the induction on
$k$. If $k=1$, then it's obvious. Assume that it's true for $k=n$,
and consider the case $k=n+1$. We can write $A$ as:
\[A=A^{\prime}\cup{A_{n+1}},\] where
$A^{\prime}=\bigcup_{i=1}^{n}{A_i}$. By the Mayer-Vietoris exact
sequence for cohomology, we have:
\[
\begin{CD}
 @. \cdots @>>> H^{d-1}(A^{\prime}\cap A_{n+1}, \mathbb{Q})\\
@>\partial>> H^d(A, \mathbb{Q}) @>i>>
H^{d}(A^{\prime},\mathbb{Q})\oplus
H^{d}(A_{n+1},\mathbb{Q}) @>j>> \cdots \\
\end{CD}
\]
It's known that the Mayer-Vietoris sequence is also an exact
sequence of mixed Hodge structures. Therefore, we get the short
exact sequence of mixed Hodge structures:
\[
\begin{CD}
0@>>>\frac{H^{d-1}(A^{\prime}\cap A_{n+1}, \mathbb{Q})}{Ker(\partial)}@>>>H^d(A, \mathbb{Q})  @>>>Ker(j) @>>>0\\
\end{CD}
\]
Notice that \[A^{\prime}\cap A_{n+1}=\bigcup_{i=1}^{n}(A_i\cap
A_{n+1}),\] therefore, by the induction assumption, both
$H^{d-1}(A^{\prime}\cap A_{n+1}, \mathbb{Q})$ and
$H^{d}(A^{\prime},\mathbb{Q})$ are Hodge-Tate structures. Now by
the Lemma \ref{h-t1}, $H^d(A, \mathbb{Q})$ is a Hodge-Tate
structure. \qed
\end{proof}

\begin{lemma}\label{h2}
Let $X$ be a complex algebraic variety, and $A$ be a subvariety of
codimension $1$. Then if both $H^{n-1}(A;\mathbb{Q})$ and
$H^n(X;\mathbb{Q})$ are Hodge-Tate, so are the relative cohomology
$H^n(X,A;\mathbb{Q})$ and $H^n(X-A;\mathbb{Q})$.
\end{lemma}
\begin{proof}
Consider the long exact sequence of cohomologies for the pair
$(X,A)$ and we know that it's also a long exact sequence of mixed
Hodge structures. Similar to Lemma \ref{h-t3},  by the Lemma
\ref{h-t1}, we can conclude that the relative cohomology
$H^n(X,A;\mathbb{Q})$ carries a Hodge-Tate structure, then by
duality, so does $H^n(X-A;\mathbb{Q})$. \qed
\end{proof}
Now let's prove the theorem:
\begin{proof}{Proof of Theorem \ref{mt}.}
By Deligne \cite{De2,De3},  there is a canonical mixed Hodge
structure on the relative cohomology
$H^{n}(\mathcal{\overline{M}}_{0,S}-A_S(\vec{a}), B_n-B_n\cap
A_S(\vec{a}))$. First we'll show that
$H^{n}(\mathcal{\overline{M}}_{0,S}-A_S(\vec{a}), B_n-B_n\cap
A_S({\vec{a}}))$ is a Hodge-Tate structure.

By the Proposition \ref{divisor-a} and Corollary \ref{nonb4}, we
know that the intersections of the components of $A_S(\vec{a})$
are either empty or Tate varieties. By Lemma \ref{h-t3} and
\ref{h2},
$H^d(\mathcal{\overline{M}}_{0,S}-A_S(\vec{a});\mathbb{Q})$
carries a Hodge-Tate structure. For the same reason,
$H^l(B_n-B_n\cap A_S(\vec{a});\mathbb{Q})$ also has a Hodge-Tate
structure. By the Lemma \ref{h2} again,
$H^{n}(\mathcal{\overline{M}}_{0,S}-A_S(\vec{a}), B_n-B_n\cap
A_S(\vec{a}))$ carries a Hodge-Tate structure. Now we need to show
that $\Omega_S{(\vec{a})}$ and $\Phi_n(\gamma)$ give the frames.
First, there are natural non-zero morphisms of pure Tate
structures by $\Omega_S{(\vec{a})}$ and $\Phi_n(\gamma)$:
\[[\Omega_S{(\vec{a})}]: \mathbb{Q}(-n)\rightarrow Gr^W_{2n}H^{n}
(\mathcal{\overline{M}}_{0,S}-A_S({\vec{a}})),\]
\[[\Phi_n(\gamma)]:\mathbb{Q}(0)\rightarrow Gr_0^{W}H_n
(\mathcal{\overline{M}}_{0,S},B_n)\]

Then composing with the canonical isomorphisms:
\[Gr^W_{2n}H^{n}(\mathcal{\overline{M}}_{0,S}-A_S({\vec{a}}),
B_n-B_n\cap A_S({\vec{a}}))\rightarrow Gr^W_{2n}H^{n}
(\mathcal{\overline{M}}_{0,S}-A_S({\vec{a}}))\] and
\[Gr_0^{W}H_{n}(\mathcal{\overline{M}}_{0,S}-A_S({\vec{a}}),
B_n-B_n\cap A_S({\vec{a}}))\rightarrow Gr_0^{W}H_n
(\mathcal{\overline{M}}_{0,S},B_n)\] we get the frame morphisms:

\[[\Omega{(\vec{a})}]^{\prime}: \mathbb{Q}(-n)\rightarrow Gr^W_{2n}
H^{n}(\mathcal{\overline{M}}_{0,S}-A_S({\vec{a}}), B_n-B_n\cap
A_S({\vec{a}}))\]

\[[\Phi_n(\gamma)]^{\prime}:\mathbb{Q}(0)\rightarrow Gr_0^{W}H_{n}
(\mathcal{\overline{M}}_{0,S}-A_S({\vec{a}}), B_n-B_n\cap
A_S({\vec{a}})).\]

Since $a_{s_i}\notin \gamma((0,1))$, for $i=1,\dots,n$, by the
description of $A_S({\vec{a}})$ and Theorem \ref{thm1} in Section
3, we have \[\Phi_n(\gamma)\cap A_S({\vec{a}})=\emptyset\]
Therefore, $\Phi_n(\gamma)$ provides an element of the relative
Betti homology
\[ [\Phi_n(\gamma)]\in H_{n}(\mathcal{\overline{M}}_{0,S}-A_S({\vec{a}}),
B_n-B_n\cap A_S({\vec{a}}))\] which is a lift of the frame
morphism $[\Phi_n(\gamma)]^{\prime}$.

The form $\Omega_S{(\vec{a})}$ obviously gives an element of the
relative De Rham cohomology class
\[[\Omega{(\vec{a})}]\in H^{n}(\mathcal{\overline{M}}_{0,S}-A_S({\vec{a}}), B_n-B_n\cap
A_S({\vec{a}})).\] And their pairing is exactly the iterated
integral in the Theorem.

Now suppose that $a_{s_i}$, $1\leq i\leq n$ are elements of a
number field $F$. We'll show that $I^{\mathcal{M}}(a_{s_1},...,
a_{s_n})$ is a framed mixed Tate motive over $F$. For this, we'll
follow the proof of \cite[Theorem~4.1]{GM} very closely. In
\cite{Gon2}, Goncharov constructed the abelian category of mixed
Tate motive over any number field $F$, where he used the theory of
triangulated category of mixed motives by Voevodsky in \cite{Vo}.
Let's apply it to our case. First, consider the standard
cosimplicial variety:
\[
S_{\bullet}(\mathcal{\overline{M}}_{0,S}-A_S({\vec{a}}),B_n-B_n\cap
A_S({\vec{a}}))\] Here
$S_0=\mathcal{\overline{M}}_{0,S}-A_S(\vec{a})$, and $S_k$ is the
disjoint union of the codimension $k$ strata of the divisor
$B_n-B_n\cap A_S({\vec{a}})$.

According to the standard procedure, we can get a complex
$S^{\bullet}(A_S({\vec{a}}),B_n)$ of varieties from the above
cosimplicial variety with $S_0$ at the degree $0$. Then it gives
an object in Voevodsky's triangulated category of mixed motives
over $F$. In fact it also belongs to the triangulated subcategory
$D_{T}(\mathbb{F})$ of mixed Tate motives over $F$. And there is a
canonical $t$-structure $t$ on $D_{T}(\mathbb{F})$. Then
$H^{n}_{t}(S^{\bullet}(A_S({\vec{a}}),B_n))$ is our mixed Tate
motive.

Similarly, by using the construction above for the Hodge-Tate
structures and the fact that the Hodge realization is a fully
faithful functor on the category of pure Tate motives, We obtain
the frame morphisms for our motive coming from
$[\Omega_S{(\vec{a})}]^{\prime}$ and $[\Phi_n(\gamma)]^{\prime}$
defined above. The theorem is proved.\qed
\end{proof}
\section{Motivic Construction of Dilogarithm}
In this section, we'll apply our preceding construction to the
dilogarithm and show that the corresponding period matrix is the
same as the one given by P. Deligne. Hence they are isomorphic.
Recall that
\[Li_{2}(z)=-\int_{0\leq t_{1}\leq t_{2}\leq
1}{\frac{dt_{1}}{t_{1}-z^{-1}}\wedge \frac{dt_{2}}{t_{2}}}\] Now
let's consider the moduli space $\mathcal{\overline{M}}_{0,5}$.
It's known that $\mathcal{\overline{M}}_{0,5}$ is the blow-up of
$\mathbb{P}^{1} \times \mathbb{P}^{1}$ at three points
$\{0,0\},\{1,1\}, \{\infty,\infty \}$ on the diagonal. Let $\pi
:\mathcal{\overline{M}}_{0,5} \rightarrow \mathbb{P}^{1} \times
\mathbb{P}^{1}$ denote this blow-up. And let $(t_{1},t_{2})$ be
the affine coordinates on $\mathbb{P}^{1} \times \mathbb{P}^{1}$.
We list the ten boundary divisors of
$\mathcal{\overline{M}}_{0,5}$ as follows (See Figure 1).
\begin{center}
\hspace{1.0cm} \epsffile{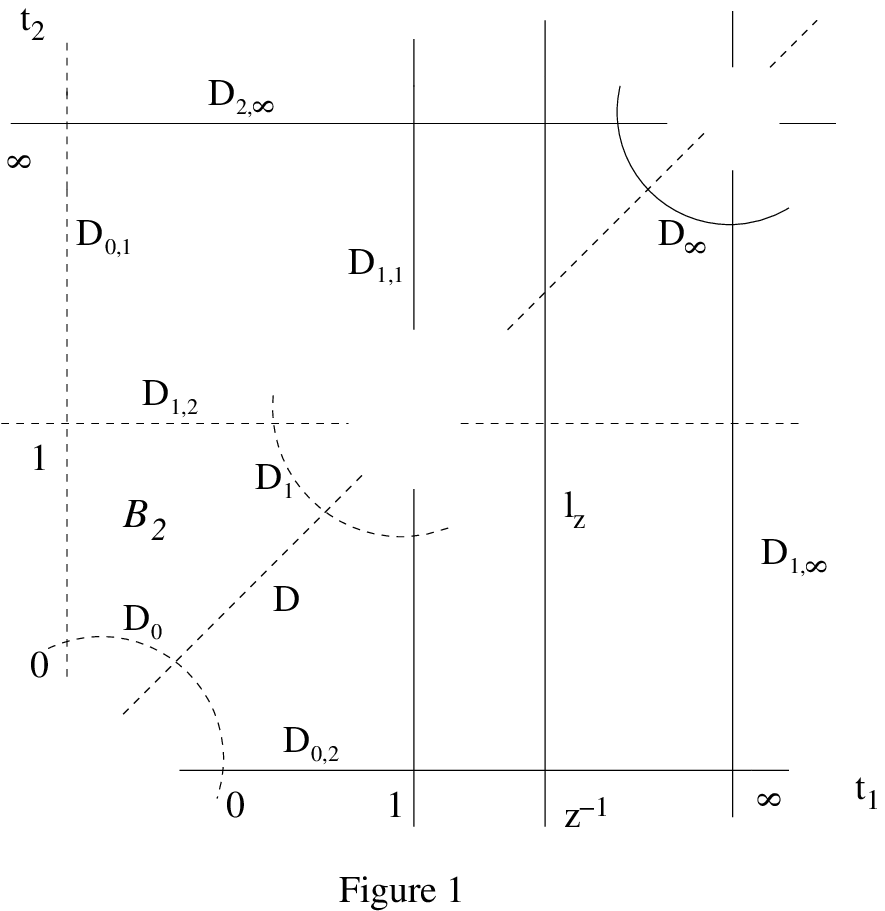}
\end{center}
\noindent $D_{1,\infty}=\text{the proper (strict) transform of the divisor}\; t_{1}=\infty ,$ \\
$D_{2,\infty}=\text{the proper transform of the divisor}\; t_{2}=\infty,$\\
$D_{\infty}=\pi^{-1}(\infty,\infty),$\\
$D_{1,1}=\text{the proper transform of the divisor}\; t_{1}=1,$ \\
$D_{1,2}=\text{the proper transform of the divisor}\; t_{2}=1,$\\
$D_{1}=\pi^{-1}(1,1),$\\
$D_{0,1}=\text{the proper transform of the divisor}\; t_{1}=0,$ \\
$D_{0,2}=\text{the proper transform of the divisor}\; t_{2}=0,$\\
$D_{0}=\pi^{-1}(0,0)$,\\
$D=\text{the proper transform of the divisor}\; t_{1}=t_{2}$.

Moreover we can describe the divisors $A(z)$ and $B_2$ explicitly.
By definition, $A(z)$ is the divisor of singularities of the
pullback meromorphic $2$-form
$\pi^{*}(\frac{dt_{1}}{t_{1}-z^{-1}}\wedge \frac{dt_{2}}{t_{2}})$
of $\mathcal{\overline{M}}_{0,5}$. We can check the following
lemma directly and leave it to the reader.
\begin{lemma}\label{l1}
\[ A(z)=
\begin{cases}
 D_{\infty,1}\cup D_{\infty,2}\cup D_{\infty}\cup D_{0,2}\cup D_{1,1},
 &\text{if $z=1$;}\\  D_{\infty,1}\cup D_{\infty,2}\cup D_{\infty}\cup D_{0,2}\cup l_{z}
 &\text{if $z\ne0,1$.}
\end{cases}
\]
where $l_{z}$ is the pullback of the divisor $t_1=z^{-1}$. That
is, $l_{z}=\pi^{-1}(t_1=z^{-1})$.
\end{lemma}

For $B_2$, by the Proposition \ref{prop1} in Section $2$, we see
that $B_2=D_{0,1}\cup D_{1,2}\cup D_{1}\cup D \cup D_{0}$ (see
Figure 1). For $z\ne0$, we'll call the following the \emph{motivic
dilogarithm}  :
\[Li_{2}^{\mathcal{M}}(z):=H^{2}(\mathcal{\overline{M}}_{0,5}-A(z),
B_2-B_2 \cap A(z)). \]
\subsection{The computations of the relative cohomology groups} First we state the results. In the following statements,
we only list the (relative) Betti homology groups because by the
comparison theorem the deRham cohomology groups are isomorphic to
the corresponding Betti cohomology groups tensor with
$\mathbb{C}$.
\begin{theorem}\label{m5.1}
 (a) If $z\ne0,1$, then the Betti homology groups of $\mathcal{\overline{M}}_{0,5}-A(z)$
are :
\[ H_{i}(\mathcal{\overline{M}}_{0,5}-A(z);\mathbb{Q})=
 \begin{cases}
  0,                            &\text{if $i\geq3$}\\
 \mathbb{Q}\oplus \mathbb{Q},  &\text{if $i=2$} \\
 \mathbb{Q},                   &\text{if $i=0,1$}

 \end{cases}
 \]
(b) If $z=1$, then the Betti homology groups of
$\mathcal{\overline{M}}_{0,5}-A(1)$ are : \[
H_{i}(\mathcal{\overline{M}}_{0,5}-A(1);\mathbb{Q})=
\begin{cases}
 \mathbb{Q},                   &\text{if $i=0,2$} \\
 0,                            &\text{otherwise}
\end{cases}
\]
\end{theorem}

\begin{theorem}\label{m5.2}
 (a) If $z\ne0,1$ then :
\[H_{2}(\mathcal{\overline{M}}_{0,5}-A(z), B_2-B_2 \cap A(z);\mathbb{Q})=
\mathbb{Q}\oplus \mathbb{Q}\oplus \mathbb{Q}\; \] Furthermore,
there are natural bases such that the corresponding period matrix
is
\[\left[
\begin{matrix}

 1         & 0               &0 \\
-Li_{1}(z) & 2\pi i          &0 \\
-Li_{2}(z) & 2\pi i\log{z}   &(2\pi i)^2

\end{matrix}
\right].\] Hence, it coincides with the one given by Deligne and
they are isomorphic as Hodge-Tate structures.

(b) if $z=1$, then:
\[H_{2}(\mathcal{\overline{M}}_{0,5}-A(1), B_2-B_2 \cap A(1);\mathbb{Q})=\mathbb{Q}\oplus\mathbb{Q}\;\]
And there are natural bases such that the corresponding period
matrix is
\[\left[
\begin{matrix}

 1          & 0       \\
 -Li_{2}(1) &(2\pi i)^2

\end{matrix}
\right].\]\\
\end{theorem}
\begin{remark}\label{rem1}
There is a dimension jump when z goes to $1$. This fact is
predicted by the specialization theorem. For more detail about it,
see \cite{Gon4}.
\end{remark}

Now we'll prove these two theorems. Before we do it, let's give a
more explicit description of the affine variety
$\overline{\mathcal{M}}_{0,5}-A(z)$. We know
$\overline{\mathcal{M}}_{0,5}$ is the blow-up of $\mathbb{P}^{1}
\times \mathbb{P}^{1}$ at the points $\{0,0\}$, $\{1,1\}$,
$\{\infty,\infty \}$, hence $\overline{\mathcal{M}}_{0,5}-(D_{1,
\infty}\cup D_{2, \infty}\cup D_{\infty})
=Bl_{\{0,1\}}(\mathbb{C}^{2})$ in which
$Bl_{\{0,1\}}(\mathbb{C}^{2})$ denotes the blow-up of
$\mathbb{C}^{2}$ at the points $\{0,0\},\{1,1\}$ (see Figure 1).
We'll use the same letter
$\pi:Bl_{\{0,1\}}(\mathbb{C}^{2})\rightarrow \mathbb{C}^{2}$ for
this blow-up.
\begin{proof}{Proof of Theorem~\ref{m5.1}:}
The idea is to decompose $\overline{\mathcal{M}}_{0,5}-A(z)$ into
the union of two subspaces, then use the Mayer-Vietoris exact
sequences.\\
\textbf{Case (a): $z\ne0,1$.} In this case we have:
\[
\overline{\mathcal{M}}_{0,5}-A(z)=(Bl_{\{0,1\}}(\mathbb{C}^{2})-\pi^{-1}\{t_{2}=0\}-\pi^{-1}\{t_{1}=z^{-1}\})\cup(D_{0}-\{*\})\notag
\]
where $*$ denotes the unique intersection point of $D_{0}$ and the
proper transform of $t_{2}=0$. Let
$\displaystyle{U_{1}=Bl_{\{0,1\}}(\mathbb{C}^{2})-\pi^{-1}\{t_{2}=0\}-\pi^{-1}\{t_{1}=z^{-1}\}}$,
then it's open and isomorphic to
$Bl_{\{1\}}(\mathbb{C}^{2}-(\{t_{2}=0\}\cup\{t_{1}=z^{-1}\}))$
which denotes the blow-up of
$(\mathbb{C}^{2}-(\{t_{2}=0\}\cup\{t_{1}=z^{-1}\})$ at the point
$(1,1)$. Let $B(\delta)=\{(t_{1},t_{2})\in \mathbb{C}^{2} : |t_1
|^{2}+|t_2 |^{2}< 4\delta^2\}$ be a small open disk around $(0,0)$
in $\mathbb{C}^{2}$. Define
$\displaystyle{U_{2}=\pi^{-1}(B(\delta))\cap(\overline{\mathcal{M}}_{0,5}-A(z))}$,
then it's clear that $D_{0}-\{*\}$ is the deformation retract of
$U_{2}$ by the line contraction: $(t_{1},t_{2})\rightarrow
(at_{1},at_{2}), 0\leq a\leq 1$. So we have the decomposition:
$\displaystyle{\overline{\mathcal{M}}_{0,5}-A(z)=U_{1}\cup
U_{2}}$. We claim that $U_{1}\cap\ U_{2}$ deformation retracts to
$S_{\delta}^{3}-S_{\delta}^{1}$, where $S_{\delta}^{k}$ is the
$k$-sphere of radius $\delta$. Indeed, notice that
$\displaystyle{U_{1}\cap\
U_{2}=\pi^{-1}(B(\delta)-\{(0,0)\})\cap(\overline{\mathcal{M}}_{0,5}-A(z))}$,
hence it is isomorphic to $B(\delta)-\{(t_{1},0):|t_{1}|< 2\delta
\}$ which deformation retracts to $S_{\delta}^{3}-S_{\delta}^{1}$
under the map $\displaystyle{f_{a}(t)=(1-a)t+ a\  \frac{\delta
\cdot t}{\|t\|},\ 0\leq a\leq 1}$, where $t=(t_1,t_2)$ and
$\|t\|=\sqrt{|t_1|^2+|t_2|^2}$. Hence,
$\displaystyle{H_i(U_{1}\cap\
U_{2})=H_i(S_{\delta}^{3}-S_{\delta}^{1})}$. By the Alexander
duality, we have
\[ H_{i}(S_{\delta}^{3}-S_{\delta}^{1};\mathbb{Q})=H^i(S^1;\mathbb{Q})=
\begin{cases}
 \mathbb{Q},   &i=0,1 \\
 0,            &\text{otherwise.}
\end{cases}
\]
The natural generator of $H_1(U_{1}\cap\ U_{2})$ is a simple loop
around the divisor $D_{0,2}$. Since $U_2$ deformation retracts to
$D_{0}-\{*\}=\mathbb{C}$, we have
$\displaystyle{H_0(U_2)=\mathbb{Q}}$ and $\displaystyle{
H_i(U_2)=0,\ i\ne0}$. Now $U_1\cong
Bl_{\{1\}}(\mathbb{C}^{2}-(\{t_{2}=0\}\cup\{t_{1}=z^{-1}\}))$, so
it follows from \cite[Page~473-474]{GH} that:
\[ H_{i}(U_{1})=
\begin{cases}
\mathbb{Q}, &i=0\\
H_1(\mathbb{C}^{2}-(\{t_{2}=0\}\cup\{t_{1}=z^{-1}\})),&i=1\\
H_{2}(\mathbb{C}^{2}-(\{t_{2}=0\}\cup\{t_{1}=z^{-1}\})\oplus
H_{2}(\mathbb{P}^{1}), & i=2\\
0, &i\geq 3
\end{cases}
\]
Now apply Mayer-Vietoris exact sequence, we obtain:\\
$\displaystyle{H_{2}(\overline{\mathcal{M}}_{0,5}-A(z))\cong H_{2}(\mathbb{C}^{2}-(\{t_{2}=0\}\cup\{t_{1}=z^{-1}\}))\oplus H_{2}(\mathbb{P}^{1})=\mathbb{Q}\oplus\mathbb{Q}\\
}$;\\
$\displaystyle{H_{0}(\overline{\mathcal{M}}_{0,5}-A(z))=\mathbb{Q};
H_{1}(\overline{\mathcal{M}}_{0,5}-A(z))=\mathbb{Q};}$
$\displaystyle{H_{i}(\overline{\mathcal{M}}_{0,5}-A(z))=0,\;
i>2.}$\\ We also see that the generator for
$H_{1}(\overline{\mathcal{M}}_{0,5}-A(z))$ is a small simple loop
around $t_1=z^{-1}$.

\textbf{Case (b):$z=1$}. We see in this case:
$\overline{\mathcal{M}}_{0,5}-A(1)=(Bl_{\{0,1\}}(\mathbb{C}^{2})-\pi^{-1}(\{t_{2}=0\})-\pi^{-1}(\{t_{1}=1\}))
\cup((D_{0}-\{*\})\cup (D_{1}-\{\bullet \}))$,
 where $*$ is the
intersection point of $D_{0}$ and $D_{0,2}$ (the proper transform
of $t_{2}=0$); $\bullet$ is the intersection point of $D_{1}$ and
$D_{1,1}$ (the proper transform of $t_{1}=1$). Let
$U_{1}=Bl_{\{0,1\}}(\mathbb{C}^{2})-\pi^{-1}(\{t_{2}=0\}\cup
\{t_{1}=1\})$. It's open and isomorphic to
$\mathbb{C}^2-(\{t_{2}=0\}\cup\{t_{1}=1\})$. As we did in Case
(a), we can construct a small open neighborhood $N_{0}$ of
$D_{0}-\{*\}$ and $N_{1}$ of $D_{1}-\{\bullet \}$ such that: (1)
$N_{0}$ deformation retracts to $D_{0}-\{*\}$, and $N_{1}$
deformation retracts to $D_{1}-\{\bullet\}$; (2) $N_{0}\cap
N_{1}=\emptyset$; (3) $N_{i}\cap U_{1}$ deformation retract to
$S_{\delta}^{3}-S_{\delta}^{1}\ (i=0,1),$ for some small $\delta$.
We define $U_{2}=N_{0}\cup N_{1}$. Then we have the decomposition:
$\displaystyle{\overline{\mathcal{M}}_{0,5}-A(1)=U_{1}\cup
U_{2}}$. And $U_{1}\cap U_{2}$ deformation retracts to the
disjoint union of two copies of $S_{\delta}^{3}-S_{\delta}^{1}$.
Now we can proceed as in Case (a) and find
$\displaystyle{H_{2}(\overline{\mathcal{M}}_{0,5}-A(1))\cong
H_{2}(\mathbb{C}^{*}\times\mathbb{C}^{*};\mathbb{Q})=\mathbb{Q}}$,
$\displaystyle{H_{0}(\overline{\mathcal{M}}_{0,5}-A(1);\mathbb{Q})=\mathbb{Q}}$,
other homology groups vanish. Theorem \ref{m5.1} is proved. \qed
\end{proof}
\begin{proof}{Proof of Theorem~\ref{m5.2}:}
Now we'll use the results of Theorem \ref{m5.1} to compute the
relative homology and cohomology groups.\\
\textbf{Case (a):$z\ne0,1$.} First, let's look at the divisor
$B_2-B_2\cap A(z)$ (See Figure 2). It consists of five irreducible
components.
Let's label them as follows (See the Figure 1 and Figure 2):\\
$l_{1}=D_{1,2}-\{z^{-1},\infty \}=\mathbb{C}-\{z^{-1}\}$,
$l_{2}=D_{1}=\mathbb{P}^{1}$, $l_{3}=D-\{z^{-1},\infty
\}=\mathbb{C}-\{z^{-1}\}$, $l_{4}=D_{0}-\{0\}=\mathbb{C}$,
$l_{5}=D_{0,1}-\{\infty \}=\mathbb{C}.$\\
We have five intersection points of these components:\\
$\displaystyle{a_{1}=l_{1}\cap l_{5};\; a_{2}=l_{1}\cap l_{2};\;
a_{3}=l_{2}\cap l_{3};\;a_{4}=l_{3}\cap l_{4};\; a_{5}=l_{4}\cap
l_{5}.\; }$
\begin{center}
\hspace{1.0cm} \epsffile{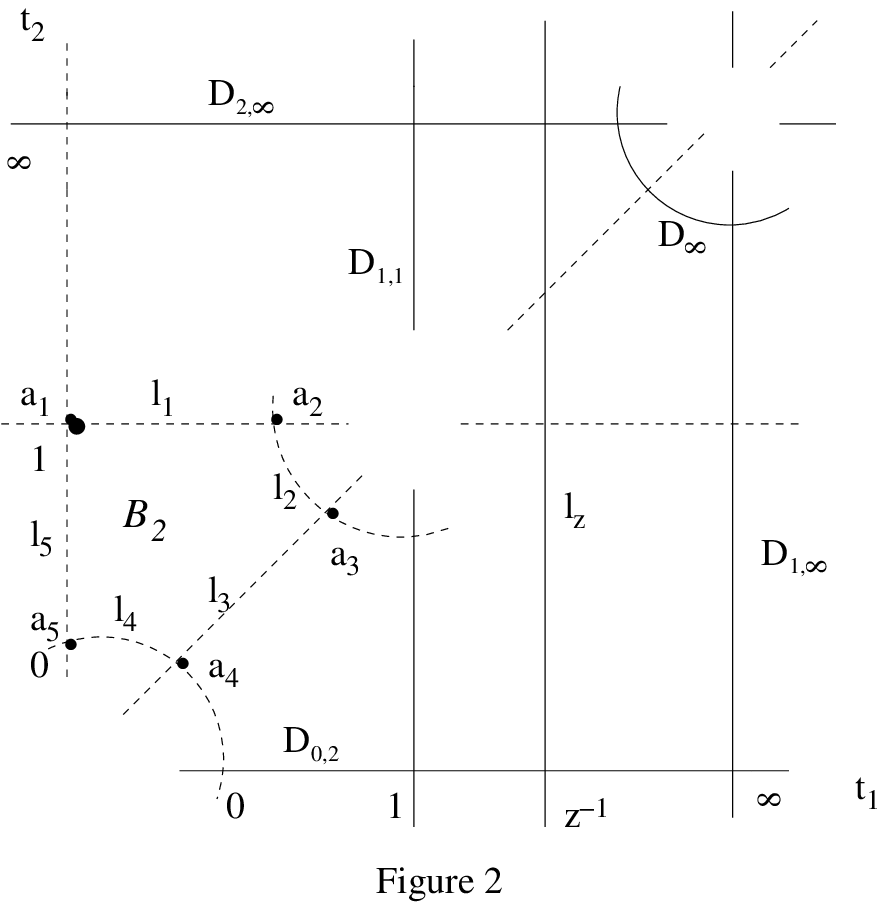}
\end{center}
Let $B_{1}=\coprod_{i=1}^5{l_i}$ be the disjoint union of the
irreducible components $l_{i}$, and $B_{0}=\coprod_{i=1}^5{a_i}$
be the disjoint union of the points $a_{i}$; $i=1,2,\dots ,5$.
Then we can compute $H_{i}(\overline{\mathcal{M}}_{0,5}-A(z),
B_2-B_2 \cap A(z))$ by the following bicomplex $(C_{p,q}, d,
\delta)$:
\[ C_{p,q}=
\begin{cases}
 C_{q}(\overline{\mathcal{M}}_{0,5}-A(z))   &\text{if $p=0$} \\
 C_{q}(B_{1})                               &\text{if $p=1$} \\
 C_{q}(B_{0})                               &\text{if $p=2$} \\
 0                                          &\text{otherwise}
\end{cases}\]
where $C_{i}(X)$ denote the vector space of singular $i$-chain on
$X$ with coefficients in $\mathbb{Q}$. The vertical differential
$d:C_{p,q}\rightarrow C_{p,q-1}$ is the differential of the
singular chain complex of $X$ and the horizontal differential
$\delta$ is defined as following:
\[\delta: C_{q}(B_{1})\rightarrow
C_{q}(\overline{\mathcal{M}}_{0,5}-A(z));\quad (c_i)_i\mapsto
\sum_{i=1}^{5}(-1)^{i-1}c_i\] where each $c_i$ is a $q$-chain of
$l_i$, $1\leq i\leq 5$. And \[\delta: C_{q}(B_{0})\rightarrow
C_{q}(B_{1});\quad \phi_{ij}\mapsto
\phi_{ij}|_{j}-\phi_{ij}|_{i}\] where for each pair $i<j$,
$\phi_{ij}$ is a $q$-chain of $l_i\cap l_j$ and $\phi_{ij}|_{i}$
represents the image of $\phi_{ij}$ in $l_i$ under the incusion
$l_i\cap l_j\hookrightarrow l_i$. That is, we have the following
diagram:
\[
\begin{CD}
\vdots                     @. \vdots                      @. \vdots \\
 @Vd VV                          @V-d VV                              @Vd VV \\
C_{2}(\overline{\mathcal{M}}_{0,5}-A(z)) @<\delta<< C_{2}(B_{1})  @<\delta<< C_{2}(B_{0})\\
 @Vd VV                                         @V-d VV          @Vd VV \\
C_{1}(\overline{\mathcal{M}}_{0,5}-A(z)) @<\delta<< C_{1}(B_{1})  @<\delta<< C_{1}(B_{0})\\
 @Vd VV                                           @V-d VV             @Vd VV \\
C_{0}(\overline{\mathcal{M}}_{0,5}-A(z)) @<\delta<< C_{0}(B_{1})
@<\delta<< C_{0}(B_{0})
\end{CD}
\]
Consider the spectral sequence with $E_{1}$ terms given as:
\[ E_{1}^{p,q}=
\begin{cases}
 H_{q}(\overline{\mathcal{M}}_{0,5}-A(z))   &\text{if $p=0$} \\
 H_{q}(B_{1})                               &\text{if $p=1$} \\
 H_{q}(B_{0})                               &\text{if $p=2$}
\end{cases}
\]
The differential $d_{1}: E_{1}^{p+1,q}\rightarrow E_{1}^{p,q}$ is
induced by the horizontal differential $\delta$. Use the explicit
generators of the homology groups in the Theorem \ref{m5.1}, we
immediately get the $E_2$ terms and $d_2=0$:
\[ E_{2}^{p,q}=
\begin{cases}
 \mathbb{Q}     &\text{if $(p,q)=(0,2),(1,1),(2,0)$} \\
 0              &\text{otherwise}
\end{cases}
\]
Hence the spectral sequence degenerates at the $E_{2}$ terms. So,
\[H_{2}(\mathcal{\overline{M}}_{0,5}-A(z), B_2-B_2 \cap
A(z));\mathbb{Q})=
\begin{cases}
\mathbb{Q}\oplus \mathbb{Q}\oplus\mathbb{Q}, &i=2\\
0, &\text{otherwise}
\end{cases}
\]
And the filtration induced by this spectral sequence coincides
with the weight filtration. Now let's turn to the relative de Rham
cohomology group $H_{dR}^{i}(\mathcal{\overline{M}}_{0,5}-A(z),
B-B \cap A(z)))$. It is defined as the $i$-th cohomology of the
total complex of the following bicomplex $(C^{p,q}, d, \delta)$:
\[ C^{p,q}=
\begin{cases}
 A^{q}(\overline{\mathcal{M}}_{0,5}-A(z))   &\text{if $p=0$} \\
 A^{q}(B_{1})                               &\text{if $p=1$} \\
 A^{q}(B_{0})                               &\text{if $p=2$} \\
 0                                          &\text{otherwise}
\end{cases}\]
where $A^{i}(X)$ denote the vector space of $C^{\infty}$  complex
$i$-forms on $X$, the vertical differential $d:C^{p,q}\rightarrow
C^{p,q+1}$ is the exterior differentiation of forms and $\delta$
is defined as
follows:\[\delta:A^{q}(\overline{\mathcal{M}}_{0,5}-A(z))\rightarrow
A^{q}(B_{1}); \quad \omega_i \mapsto (-1)^{i-1}\omega|_{l_i}\]
where $\omega$ is a $q$-form of
$\overline{\mathcal{M}}_{0,5}-A(z)$ and $\omega|_{l_i}$ is the
restriction of $\omega$ on $l_i$.
\[\delta: A^{q}(B_{1})\rightarrow A^{q}(B_{0}); \quad \delta\theta \mapsto
\theta_j|_{ij}-\theta_i|_{ij}\] where $\theta=(\theta_i)$ is a
$q$-form of $B_1$ and for each pair $i<j$, $\theta_i|_{ij}$ is the
restriction of $\theta_i$ on $l_i\cap l_j$. We have the following
diagram:
\[
\begin{CD}
\vdots                     @. \vdots                      @. \vdots \\
 @AdAA                          @A-dAA                              @AdAA \\
A^{2}(\overline{\mathcal{M}}_{0,5}-A(z)) @>\delta>> A^{2}(B_{1}) @>\delta>> A^{2}(B_{0})\\
@AdAA                                              @A-dAA             @AdAA \\
A^{1}(\overline{\mathcal{M}}_{0,5}-A(z)) @>\delta>> A^{1}(B_{1}) @>\delta>> A^{1}(B_{0})\\
@AdAA                                           @A-dAA               @AdAA \\
A^{0}(\overline{\mathcal{M}}_{0,5}-A(z)) @>\delta>> A^{0}(B_{1})
@>\delta>> A^{0}(B_{0})
\end{CD}
\]

Similar to the case for homology, we see that the corresponding
spectral sequence degenerates at the $E_{2}$ terms. So we obtain:
$H_{dR}^{2}(\mathcal{\overline{M}}_{0,5}-A(z), B-B \cap
A(z)))=\mathbb{C}\oplus \mathbb{C}\oplus \mathbb{C}$, and other
relative cohomology groups vanish. For the case $z=1$, the proof
is similar and we leave it to the reader.
\subsection{Computation of the Period Matrix}
We have the following three cochains in the total complex
associated to the deRham bicomplex $(C^{p,q}, d, \delta)$:
\[e_{1}=(0,0,\delta_{A_{1}}); e_{2}=(0,
\frac{dt_{1}}{t_1-z^{-1}},0);
e_{3}=(\pi^{*}(\frac{dt_{1}}{t_1-z^{-1}}\wedge
\frac{dt_{2}}{t_{2}}),0,0)\] where $\delta_{A_{1}}$ is a function
on $B_{0}$ satisfying
$\delta_{A_{1}}(A_{1})=1,\delta_{A_{1}}(A_{j})=0, j\ne1$ and
$\frac{dt_{1}}{t_1-z^{-1}}$ is a $1$-form on the component $l_1$
of $B_1$. Since $\frac{dt_{1}}{t_1-z^{-1}}$ and
$\pi^{*}(\frac{dt_{1}}{t_1-z^{-1}}\wedge \frac{dt_{2}}{t_{2}})$
are holomorphic on $l_1$ and $\mathcal{\overline{M}}_{0,5}-A(z)$
respectively, they are closed forms. Hence $e_{i}$ are cocycles
and represent elements in
$H_{dR}^{2}(\mathcal{\overline{M}}_{0,5}-A(z), B-B \cap A(z))$.\\
Next, let's consider the cycles in the Betti bicomplex $(C_{p,q},
d, \delta)$.\\ Let
$\displaystyle{b_1=(\overline{\Phi}_{2},\partial{(\overline{\Phi}_{2})},\sum_{i=1}^{5}{a_i})}$,
where $\partial{(\overline{\Phi}_{2})}$ denotes the boundary of
$\overline{\Phi}_{2}$ in $B_1$. It's just the $5$ sides of the
dotted pentagon in figure $2$. Let
$\displaystyle{b_2=(C,\partial{C},0)}$, where $C=\pi^{-1}(C_0),
C_0=(z^{-1}+\epsilon e^{2\pi iu}, v+(1-v)(z^{-1}+\epsilon e^{2\pi
iu}))\subset \mathbb{C}\times\mathbb{C}\subset \mathbb{P}^1 \times
\mathbb{P}^1$, and $\epsilon$ is a small positive number, $0\leq u
\leq 1, 0\leq v \leq 1$. It's clear that the boundary
$\partial{C_0}$ consists of $2$ small cycles $(z^{-1}+\epsilon
e^{2\pi iu},1)$ and $(z^{-1}+\epsilon e^{2\pi iu},z^{-1}+\epsilon
e^{2\pi iu})$. Since $C_0$ does not contain the blow-up points,
$C$ is isomorphic to $C_0$ and the boundary
$\partial{C}=\partial{C_0}$. The key point is that $\partial{C}$
is contained in $B_{1}=\coprod_{i=1}^5{l_i}$, which means exactly
that $b_2$ is a cycle in the total complex. Finally, let
$\displaystyle{b_3=(T,0,0)}$, where $T=\pi^{-1}(z^{-1}+\epsilon
e^{2\pi iu},\epsilon e^{2\pi iv}), (0\leq u \leq 1, 0\leq v \leq
1)$ is the inverse image of a torus and since the three blow-up
points are not on the torus, $T$ is isomorphic to its image.
Clearly the boundary of $T$ is zero. Hence $b_3$ is a cycle.

Now we can calculate the period matrix $P=(p_{ij})$ between
$(b_1,b_2,b_3)$ and $(e_1,e_2,e_3)$, here $p_{ij}=<e_i,b_j>$, $1\leq
i, j\leq3$. It's straightforward that $<e_1,b_1>=1$,
$\displaystyle{<e_2,b_1>=-Li_{1}(z),}$ \\$\displaystyle{
<e_3,b_1>=-Li_{2}(z)}$, and $<e_1,b_3>=0$,
$\displaystyle{<e_2,b_3>=0, <e_3,b_3>=(2\pi i)^2, <e_1,b_2>=0}$.
Since $<e_2,b_2>$ is equal to the integral of the $1$-form
$\displaystyle{\frac{dt_{1}}{t_1-z^{-1}}}$ over a small circle with
center $t_1=z^{-1}$, we have $\displaystyle{<e_2,b_2>=2\pi i}$. And
\begin{alignat}{2}
<e_3,b_2>&= \int_{C}\pi^{*}(\frac{dt_{1}}{t_1-z^{-1}}\wedge
            \frac{dt_{2}}{t_{2}})=\int_{C_0}\frac{dt_{1}}{t_1-z^{-1}}\wedge
            \frac{dt_{2}}{t_{2}}\notag\\
         &=\int_{0}^{1} \int_{0}^{1}2\pi i du \cdot
           d(\log{(v+(1-v)(z^{-1}+\epsilon e^{2\pi iu})}))\notag\\
         &=-2\pi i\int_{0}^{1}\log{(z^{-1}+\epsilon e^{2\pi iu})}du=-2\pi i\log{z^{-1}}&& 
         \notag\\
         &=2\pi i\log{z}\notag
\end{alignat}

Therefore, we get the period matrix:
\[\left[
\begin{matrix}

 1         & 0               &0 \\
-Li_{1}(z) & 2\pi i          &0 \\
-Li_{2}(z) & 2\pi i\log{z}   &(2\pi i)^2
\end{matrix}
\right].\] Since it's nonsingular, it follows that $\{e_1, e_2,
e_3\}$ and $\{b_1, b_2, b_3\}$ are bases. For part(b), it's clear
that: for the cocycles, we have $e_{1}=(0,0,\delta_{A_{1}}),
e_{3}=(\pi^{*}(\frac{dt_{1}}{t_1-z^{-1}}\wedge
\frac{dt_{2}}{t_{2}}),0,0)$; for the cycles, we have
$b_1=(\overline{\Phi}_{2},\partial{(\overline{\Phi}_{2})},\sum_{i=1}^{5}{a_i}),
b_3=(T,0,0)$. Therefore, its period matrix is\\
\[\left[
\begin{matrix}

 1          & 0       \\
 -Li_{2}(1) &(2\pi i)^2
\end{matrix}
\right].\] \qed
\end{proof}

\end{document}